\documentclass[]{tMAM2e}

\usepackage{epstopdf}
\usepackage{graphicx}
\usepackage{epstopdf}
\usepackage{psfrag}
\usepackage{hyperref}
\usepackage[section]{placeins}
\hypersetup{
    colorlinks,
    citecolor=black,
    filecolor=black,
    linkcolor=black,
    urlcolor=blue
}

\def\Real{\mathbb{R}}
\def\Zed{\mathbb{Z}}

\def\kemal{Karaosmano\u{g}lu}

\begin{document}

\doi{10.1080/1745973YYxxxxxxxx}

\jnum{Number} \jyear{2013} \jmonth{Month}

\markboth{W. A. Sethares and R. Budney}{Journal of Mathematics and Music}

\articletype {RESEARCH ARTICLE}

\title{Topology of Musical Data}

\author{William A. Sethares$^{\rm a}$
\vspace{6pt} and Ryan Budney$^{\rm b}$\\
\vspace{9pt}  $^{\rm a}${\em{Department of Electrical and Computer Engineering, University of Wisconsin, Madison, USA, sethares@wisc.edu}}; $^{\rm b}${\em{Mathematics and Statistics, University of Victoria, PO Box 3060 STN CSC, Victoria, B.C., Canada V8W 3R4, rybu@uvic.ca}}\\
\received{4 Jan. 2013; final version received d mmmm yyyy} }

\maketitle

\begin{abstract}
The musical realm is a promising area
in which to expect to find nontrivial topological structures. 
This paper describes several kinds of metrics on musical data,
and explores the implications of these metrics in two ways: via techniques
of classical topology where the metric space of all-possible musical 
data can be described explicitly, and via modern data-driven
ideas of persistent homology which calculates the Betti-number bar-codes
of individual musical works.
Both analyses are able to recover three well known topological structures 
in music: the circularity of octave-reduced musical scales, the
circle of fifths, and the rhythmic repetition of timelines.
Applications to a variety of musical works 
(for example, folk music in the form of standard MIDI files) are
presented, and the bar codes show many interesting features.
Examples show that individual pieces may span the complete space (in
which case the classical and the data-driven analyses agree), or they may span 
only part of the space.
\end{abstract}

\begin{keywords}persistent homology, 
topological structure,
metrics on musical data,
melodic analysis,
rhythm analysis,
circle of fifths
\end{keywords}

\maketitle

\section{Introduction}

Music has rich internal structures.
Embedded in any visualization of musical structure 
is a notion of the closeness and/or similarity (or equivalently
a notion of the distance and/or dissimilarity) between various musical elements.
This paper considers the implications of
several different distance-measures applied at the level of notes,
of chords/scales, and of rhythms. Two themes emerge.
Ideas of musical similarity tend to be compound,
to express relationships more complex than are easily visualized in
(linear) Euclidean space; each notion of closeness corresponds to
a space with different topological structure. 
Second, even (apparently) small changes in the definition of 
a distance function can lead to large changes in the global
(topological) shape of the resulting space. 
The implication is that there is not really a single ``master''
space in which all musical phenomena can be placed. Rather,
different notions of distance between different kinds of 
musical events can be legitimately described in various ways,
as having their own kind of local and global structure.

A large literature shows many ways to visualize and analyze musical structures.
To name a few examples: musical scores present notes in a two-dimensional array,
Euler's classical {\em Tonnetz} represents musical intervals in a lattice as do modern generalizations
\cite{catanzaro}, Partch \cite{partch} draws the tonality diamond where the two axes 
represent significant musical intervals, 
Chew \cite{chew} visualizes musical progressions along a spiral array, and
computer-based music visualizers \cite{musicVisualization}
display moving patterns in real time as
music progresses. Lewin \cite{lewin} suggests the application of simplicial complexes and
homology to music theory and Mazzola \cite{mazzola1990} considers
actions by the symmetric group on $n$-tuples. More recently,
Callender, Quinn, and Tymoczko
\cite{callender2005}-\cite{callender2008} explore a topological approach to musical
spaces using the voice-leading metric; the result is a space they call
an orbifold. 
Topological considerations have also been used 
in \cite{buteau2008} and \cite{buteau2004} to cluster
melodies and so play a role in the discovery of musical motives.

Several different notions of distance applicable to musical situations 
are discussed in Section \ref{sec:distances} and are
intended to motivate both the mathematical analyses in Section \ref{sec:symm}.
Section \ref{sec:persHomology} reviews the idea of persistent homology,
which is then applied to musical data in Section \ref{sec:simulations}.
 
Throughout the paper, several interesting 
spaces are encountered, not all of which can be described as open subsets of $\Real^{n}$. 
Each of the distance functions corresponds to a particular
topological space, and these spaces can often be described concretely
as products of basic spaces like the circle $S^{1}$ and the closed ball $D^{n}$. 
Occasionally, somewhat more elaborate spaces arise, such as twisted products
(fibre bundles) and higher-dimensional spheres.

A collection of simulations examine a series of increasingly complex sets of musical data using 
the Betti numbers of persistent homology. Applications begin with
simple note and chord sequences to verify that the method is capable of displaying common-sense
musical intuition. In some of the examples, the song data 
spans much of the underlying topological space while in 
other examples the musical data occupies only a small part of the 
space. Both situations may be interesting.
Finally, an online supplement \cite{onlineSupplement} is available
which contains additional figures and examples that augment the
presentation within the paper.

\section{Distance Measures} \label{sec:distances}

\begin{figure}
\begin{center}
\includegraphics[width=1.8in]{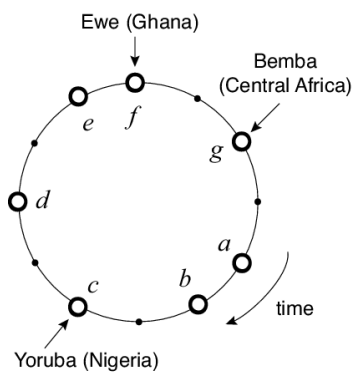}
\caption{In the necklace notation, time progresses around the circle in a 
clockwise direction, marked off by a series of circles: small
dots represent silent time points while open circles represent
sounding events. Each passage around the circle represents one cycle through
the rhythm. Traditional rhythms of the Ewe (from Ghana), the Yoruba (from
Nigeria), and the Bemba (from Central Africa) are all variants
of the ``standard rhythm pattern'' described by King \cite{king61}, and
can be represented in the necklace notation by the same pattern
but with different starting points.}
\label{fig:yoruba}
\end{center}
\end{figure}

One of the oldest visualizations of musical structure
represents temporal cycles as spatial circles:  Saf\^i al-Din al-Urmaw\^i,
the 13th century theoretician from Baghdad,
represents both musical and natural rhythms in a circular notation 
in the {\em Book of Cycles} \cite{al-din}.
Time moves around the circle (usually in 
a clockwise direction) and events are depicted along the periphery.
Since the ``end'' of the circle is also the ``beginning,'' this
emphasizes the repetition inherent in rhythmic patterns. 
Anku \cite{anku} argues that African music is perceived in a 
circular (rather than linear) fashion that makes the necklace notation,
shown here in Figure \ref{fig:yoruba}, particularly appropriate.

In this setting, musical events are sounds occurring at times
specified by their position around the circumference of the circle. 
There are two different notions of the closeness of events in this
representation. In one notion, two events are close
if they occur near each other in time (for instance, events represented by the points labelled
$a$ and $b$ are closer together than events represented by the points $c$ and $d$.
In the second notion, events are close if they occur at (or near) the
same position in each cycle. Thus an event at $a$ (which occurs the first time through the cycle) is identical to another event that occurs in a later cycle at $a$. With only the first notion
of proximity in time, the natural representation would be a line segment $I^{1}=[0,1]$;
with the repetitions of the cyclical pattern, the line 
is transformed into a circle $S^{1}$.

To be concrete, the distance between two events occurring at 
times $f$ and $g$ is
\begin{equation} \label{eqn:distanceRhy}
d_{N}(f,g) = \min(s,1-s) \mbox{ where } s= |f-g| \hspace{-3mm} \mod 1,
\end{equation}
and where one unit of time represents one period of the rhythm.
Call this the {\em necklace distance}, indicated by the subscript $N$. 
The necklace metric on the space of repetitive rhythms is equivalent to the natural metric
on the quotient space $S^{1} = \Real / \Zed$. The slash notation indicates the conversion of the
linear structure of the reals to a cyclic structure, by considering two numbers equivalent if they differ by an integer. 

\begin{figure}
\begin{center} 
\includegraphics[width=4.0in]{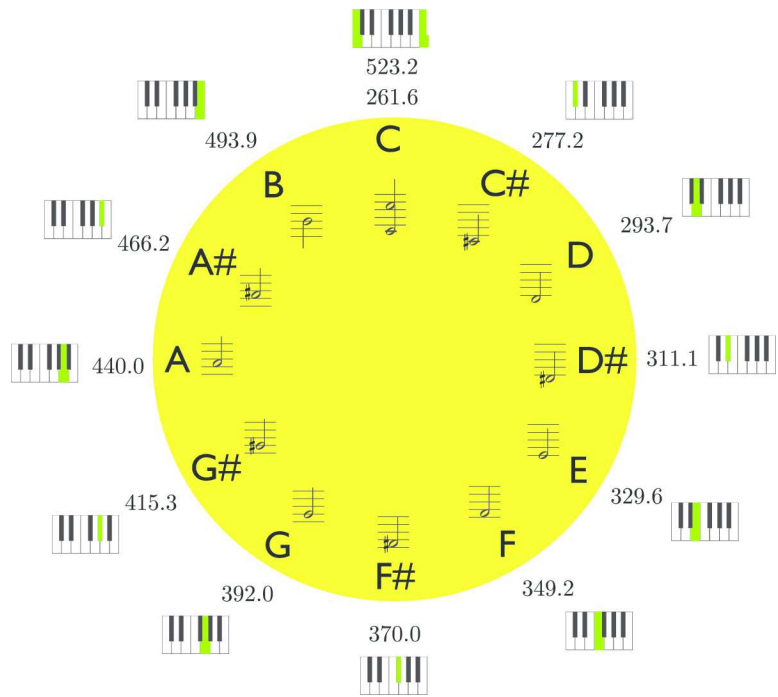}
\end{center}
\caption{One way to picture
two important aspects of musical perception is to place notes that are near each other in frequency
(such as $C$ and $C \sharp$) around the edge of a circle, and to overlay 
notes that are an octave apart (i.e, those which are distance zero 
apart in the pitch class metric (\ref{eqn:distancePC}))
such as low $C$ and high $C$.}
\label{fig:circle-of-notes}
\end{figure}

Music, of course, has aspects other than rhythm.
The fundamentals of complex tones are generally perceived as a function 
of frequency in a logarithmic fashion. Thus the $15.6$ Hz ``distance'' from 
$C$ to $C\sharp$ is perceived to be the same size as the $24.7$ Hz ``distance''
from $G\sharp$ to $A$ (refer to Figure \ref{fig:circle-of-notes} for
the origin of these numerical values). Accordingly, it is reasonable to
consider a measure that operates on $\log$ frequency rather than
on frequency itself. A metric like $|\log_{2}(f)-\log_{2}(g)|$
captures this.
But there is also a second notion of closeness operating in 
the pitch domain, on the set of positive reals $\mathbb{R}_+$, 
that of octave equivalence. For example,
the $C$ note at $261.6$ Hz and the high $C$ at $523.2$ Hz are 
closely identified. This is what happens when a man sings 
along with a woman (or when a woman sings along with a child):
the ``same'' note is actually a factor of two apart in frequency.
These two notions of nearness combine to suggest measuring distance
between frequencies as
\begin{equation} \label{eqn:distancePC}
d_{PC}(f,g) = \min(s,1-s)
\mbox{ where } s = |\log_{2}(f)-\log_{2}(g)| \hspace{-3mm}\mod 1.
\end{equation}
For example, the distance $d_{PC}(261.6 Hz, 523.2 Hz)$ 
between the two C notes an octave apart, is zero. This means that $d_{PC}$ is not a metric 
(because metrics must satisfy the ``identity of indiscernibles,'' that is, 
$d(x,y)=0 \Leftrightarrow x=y$).
But (\ref{eqn:distancePC}) does induce a metric on the quotient space determined by the 
equivalence relation $x \sim y \Leftrightarrow d_{PC}(x,y)=0$.
In musical terms, this may be interpreted as a measure of distance between ``pitch classes''
\cite{rahn}, since it identifies all $C$s, all $C\sharp$s, etc. into 
equivalence classes and so (\ref{eqn:distancePC}) is called 
the {\em pitch-class} distance on $\mathbb{R}_+$, as indicated by the
subscript $PC$. Again, two notions of
closeness have combined to turn the space of pitch-class notes into
the natural metric on the quotient space $S^{1} = \Real / \Zed$.

Of course, music contains more than just isolated pitch classes and 
temporal cycles. One way to combine time and pitch information 
is with a time-delay embedding, which is a common practice in time series analysis.
Suppose that a melody consists of a sequence of notes with fundamental frequencies 
at $f_{1}, f_{2}, f_{3}, f_{4}... $.
These may be combined into pairs (a two-dimensional time-delay
embedding) by creating the sequence of vectors $(f_{0}, f_{1})$, 
$(f_{1}, f_{2})$, $(f_{2}, f_{3})$, $(f_{3}, f_{4})...$.
The distances between such vectors can be calculated by 
adding the distances between the pitches element-wise 
using the pitch-class metric (\ref{eqn:distancePC}). 
For example, the distance between  
pairs $(f_{n}, f_{n+1})$ and $(f_{m}, f_{m+1})$ 
is $d_{PC}(f_{n},f_{m})+d_{PC}(f_{n+1},f_{m+1})$.  
The collection of all pairs of frequencies under this distance measure can 
be thought of as the torus $S^1 \times S^1$.
Similarly, a $d$-dimensional time-delay embedding 
(where the $f_{i}$ are gathered into vectors of size $d$)
can be visualized as the $d$-fold product 
\[
\overset{d \ \small{ times}}{\overbrace{S^1 \times S^1 \times ...  \times S^1 \times S^1}}.
\]

It may also be useful to generalize the distance measure to
consider scalar harmony consisting of multiple pitches sounding simultaneously.
The most straightforward generalization of the pitch-class
metric is to add the (pitch-class) distances between all 
elements of the vectors, as was done for the time-delay
embeddings. This distance measure distinguishes
chord inversions: for instance, a $C$-major chord in root position
($C$-$E$-$G$) would be distant from a $C$-major chord in 
third position ($G$-$C$-$E$). While this may be desirable
in some musical situations, it is undesirable when looking for
structures that involve musical key, where (say) all $C$-major chords 
should be identified irrespective of inversion and all $C$ major scales are 
identified irrespective of the order in which the pitches are listed.
For example, the ascending $C$-major scale 
and the descending $C$-major scale are both the same
entity, and the distance should reflect this
realm of musical perception. 
Another example showing the shortcomings of (\ref{eqn:distancePC})
is when two voices are sounding simultaneously.
The pairs $(f_{1}, f_{2})$ and $(f_{2}, f_{1})$
both mean that same two notes are heard. Using the 
pitch-class metric, these are unequal.

This can be addressed by permuting the elements and calculating 
the minimum over all the permutations.
To be precise, let $f=(f_{1}, f_{2}, ... f_{n})$ and 
$g=(g_{1}, g_{2}, ... g_{n})$ be two $n$-tuples, and 
define the distance
\begin{equation} \label{eqn:distanceCC}
d_{CC}(f,g) = \min_{\sigma} \sum_{i=1}^{n} d_{PC}(f_{i}, \sigma (g_{i})) 
\end{equation}
where $\sigma$ ranges over all possible permutations of size $n$, and 
$d_{PC}(\cdot, \cdot)$ is the pitch-class distance of 
(\ref{eqn:distancePC}).
This {\em chord-class distance} is the smallest of the pitch-class
distances between $f$ and all permutations of $g$. Hence it
is invariant with respect to chord and scale inversion;
all reorderings of the elements of $f$ and $g$ are
placed in the same equivalence class (and have zero distance
from each other). 
The quotient space induced by the chord-class distance is called the
{\it symmetric product}.  Specifically, given a space $X$, the $n$-fold 
symmetric product is denoted $Symm_n(X) = X^n / \Sigma_n$, and $\Sigma_n$ is
the symmetric group acting on $X^n$ by permuting the coordinates of $X^n$. 
The chord-class distance induces a metric on $Symm_n(S^1)$, as described in
Section \ref{sec:symm}.  
A further extension of these ideas combines a length $j$ time-delay embedding of $k$-note chords.
Describing the structure of these spaces concretely requires some work as in
Section \ref{sec:symm}.
Another complication is that there may not always be exactly 
$n$ elements in each of the chords: for instance, one of the voices may
become momentarily silent or two of the voices may converge to the same pitch.
This can be handled by exploiting the $n$-th finite subset subspace, 
and leads to some of the more complex spaces encountered in the analysis.

As a further extension to rhythmic patterns such as those of Figure \ref{fig:yoruba},
consider extending the distance metric to one which measures the
distance between complete rhythmic sequences (rather than
the distance between individual timed events). In this case, an appropriate measure
is one which considers any two sequences equal if they are circular shifts
of each other. Let $f$ and $g$ be $n$-vectors that define  
the times of events in a cycle with period 1.
Then the {\em rhythm distance} between $f$ and $g$ is
\begin{equation} \label{eqn:distanceCyc}
d_{R}(f,g) = \min_{\sigma_{c}} \sum_{i=1}^{n} d_{N}(f,\sigma_{c}(g)),
\end{equation}
where $\sigma_{c}$ ranges over all cyclic shifts and $d_{N}(\cdot,\cdot)$ is the 
distance defined in (\ref{eqn:distanceRhy}). Under this measure, all three of the 
named rhythms in Figure \ref{fig:yoruba} are the same since they are cyclical shifts
of each other. In this case, the topology will be given by the subsymmetric product 
$Symm_n^G(X) = X^n / G$, which is the quotient of the $n$-times product of the 
circle $X=S^{1}$ and the subgroup $G \subset \Sigma_n$ of circular shifts.

Finally, a standard metric that can be applied to finite subset spaces is the {\it Hausdorff distance} 
which defines the distance between two finite subsets of a metric space $A,B \subset X$ as
 
\[
d_H(A,B) = \max \{ \max_{a \in A} \{ \min_{b \in B} d(a,b) \} \}, \max_{b \in B} \{ \min_{a \in A} d(a,b) \} \}.
\]

\section{Spaces and Homology}
\label{sec:symm}

This section provides the tools needed to interpret the examples in Section
\ref{sec:simulations}, which analyze musical data using persistent homology.  
Sections \ref{sec:def}, \ref{sec:spaces}, and \ref{sec:finiteSubspaces} describe auxiliary {\it natural}
spaces on which it is possible to analyze the input (musical) data.  These spaces admit their own notion of (metric) geometry compatible with Section \ref{sec:distances}, and it is through the lens of these geometries that persistent homology will be applied.  The use of such auxiliary geometries is a commonly-occurring theme in topology. For example, \cite{CZCG} gives a detailed discussion and connections to subjects such as the {\it Shape Theory} of \cite{kendall}. Section \ref{sec:homol} summarizes the basics of homology in and Section \ref{sec:persHomology} considers persistent homology.  The examples of Section \ref{sec:simulations} show that the persistent homology of data is influenced by the topological features of the ambient metric spaces in which they live. 

\subsection{Products, Symmetric Products and Bundles} \label{sec:def}

Given a pair of topological spaces $X$ and $Y$, the {\it product} of the two
is denoted $X \times Y = \{ (x,y) : x \in X, \text{ and } y \in Y \}$.  For 
example, if $S^1$ is the circle, $S^1 \times S^1$ is called a {\it torus}.  
Similarly, $S^1 \times \{1,2\}$ is {\it two circles}.  The $n$-fold product of $X$, denoted $X^n$, 
is the set of all $n$-tuples of elements of $X$. Thus 
$X^n = \{ (x_1,\cdots,x_n) : x_i \in X \ \forall i \in \{1,2,\cdots ,n\}\}$.  
Observe that elements of $X^n$ consist of $n$ ordered elements of $X$.   
For example, $(x_1,x_2,\cdots,x_n) = (y_1,y_2,\cdots, y_n)$ if and only if 
$x_1=y_1, x_2=y_2, \cdots, x_{n-1}=y_{n-1}$ and $x_n=y_n$. 

The $n$-fold {\it symmetric product} of $X$ is the quotient of $X^n$, 
by the action of the {\it symmetric group} $\Sigma_n$, that is, $Symm_n(X) = X^n / \Sigma_n$ where  
$\Sigma_n \equiv \Sigma \{1,2,\cdots,n\}$ is the set of all one-to-one onto 
functions $\{1,2,\cdots,n\} \to \{1,2,\cdots,n\}$ acting on $X^n$ on the right.  
Elements of $Symm_n(X)$ consist of $n$ unordered elements of $X$ and so  
$(x_1,x_2, \cdots, x_n) = (y_1,y_2,\cdots,y_n)$ in $Symm_n(X)$ if and only if there is a 
permutation $\sigma \in \Sigma_n$ such that $x_i = y_{\sigma(i)}$ for all $i \in \{1,2,\cdots,n\}$. 
Given an arbitrary subgroup $G \subset \Sigma_n$, the {\it sub-symmetric product} $Symm_n^G(X)$ is the quotient $Symm_n^G(X) = X^n / G$.   In $Symm_n^G(X)$, $(x_1,x_2,\cdots, x_n) = (y_1,y_2,\cdots, y_n)$ if and only if there is a 
permutation $\sigma \in G$ such that $x_i = y_{\sigma(i)}$ for all $i \in \{1,2,\cdots, n\}$. 

A frequently-useful variant on the idea of the product of two spaces is the notion of a {\it fibre bundle}. A fibre-bundle is a continuous function $p : E \to B$ which is `locally trivial' in the sense that for every $x \in B$ there is a neighbourhood $U \subset B$ and a homeomorphism $\phi : p^{-1}(U) \to U \times F$ where $F = p^{-1}(x)$.  $F$ is called `the fibre over $x$' of the fibre bundle.  The map $\phi$ must satisfy that $\pi \circ \phi = p$, where $\pi : U \times F \to U$ is the projection map onto the first coordinate $\pi(u,f) = u$.  

A product $X \times Y$ is a fibre bundle with $E = X \times Y$, $B = X$ and $F = Y$, with $p : X \times Y \to X$ given by $p(x,y) = x$.  One should think of a fibre bundle as a space that is locally a product in that the map $\phi$ gives a way of locally identifying nearby fibres.  The M\"obius band is an example of one such space, as we will see in Section \ref{sec:spaces}.  A fibre bundle with fibre a disc is called a {\it disc bundle}. 

\subsection{Morton's Bundles} \label{sec:spaces}

A classical result of Morton \cite{Morton} states that $Symm_n(S^1)$ is a disc bundle over $S^1$ which is orientable if and only if $n$ is odd.  This means that when $n$ is odd, the local homeomorphisms $\phi$ from the definition of a bundle assemble to give a homeomorphism $Symm_n(S^1) \cong S^1 \times D^{n-1}$.  The space $Symm_2(S^1)$ is a M\"obius band, and when $n$ is even, $Symm_n(S^1)$ is similar to a M\"obius band in that the local homeomorphisms of the fibres assemble into a mirror reflection of the fibre.  

Morton gives a general prescription for studying sub-symmetric products $Symm_n^G(S^1)$.  The basic idea is to consider the equivalence of $S^1$ and $\Real / \Zed$.  This makes it possible to think of $(S^1)^n / G$ as $\Real^n / (\Zed^n \rtimes G)$ where $\Zed^n \rtimes G$ denotes the semi-direct product\footnote{Given a group $G$ and a homomorphism $\phi : G \to Aut(H)$ the semi-direct product $H \rtimes_\phi G$ is defined as pairs $(h,g)$ with multiplication defined as $(h_1,g_1)\cdot(h_2,g_2) = (h_1\phi_{g_1}(h_2), g_1g_2)$.} of $G$ with $\Zed^n$, equivalently it is denoted $\Zed \wr G$ and denotes the wreath product\footnote{The wreath product $\Zed \wr G$ is a synonym for $\Zed^n \rtimes G$ where $G$ acts on $\Zed^n$ via the regular representation $G \subset \Sigma_n$.} of $\Zed$ with $G$.  Complex multiplication gives a $G$-equivariant map $(S^1)^n \to S^1$ where the action of $G$ on the target $S^1$ is trivial.   Similarly, $G$ acts on $\Real^n$ via the inclusion $G \hookrightarrow \Sigma_n$ making $\Real^n \to (\Real / \Zed)^n$ $G$-equivariant.  
Denoting the kernel of the addition function 
$+ : \Real^n \to \Real$ by $\Delta_\Real = \{(x_1,x_2,\cdots,x_n) \in \Real^n : x_1+x_2+\cdots+x_n=0\}$ 
and $\Delta_\Zed = \Zed^n \cap \Delta_\Real$, this means that there is a fibre bundle
\[
\Delta_\Real / (\Delta_\Zed \rtimes G) \to \Real^n / (\Zed^n \rtimes G) \to \Real/\Zed .
\]
Moreover, Morton identifies the monodromy of the bundle.

\begin{theorem}\label{MortonThm} 
The space $Symm_n^G(S^1)$ is a fibre bundle over $S^1$, with fibre $\Delta_\Real / (\Delta_\Zed \rtimes G)$ and monodromy induced by the map $(t_1,t_2, \cdots, t_n) \longmapsto (t_1 + \frac{n-1}{n}, t_2 - \frac{1}{n}, \cdots, t_n -\frac{1}{n})$ of $\Delta_\Real$. 
\end{theorem}

The other main result in Morton's paper is the identification of the quotient $\Delta_\Real / (\Delta_\Zed \rtimes G)$ as a disc in the case $G =\Sigma_n$.  A key argument is that the region $\{ (t_1,\cdots,t_n) \in \Delta_\Real : t_1 \leq t_2 \leq \cdots \leq t_n \leq t_1+1\}$ is a fundamental domain for the action of $\Delta_\Zed \rtimes \Sigma_n$ on $\Delta_\Real$.  Since there are no identifications on the boundary, this {\it is} the fibre.  It's an exercise to check the monodromy is orientation preserving if and only if $n$ is odd.  There are only two disc bundles over a circle: the trivial product bundle $D^{n-1} \times S^1$ and the M\"obius bundle $D^{n-1} \rtimes S^1$.    The `center' ($0$-section) of the disc bundle consists of the equally-spaced points $(t_1, t_1 + \frac{1}{n}, t_1+\frac{2}{n}, \cdots, t_1 + \frac{n-1}{n})$. 

Depending on the group $G$, $Symm_n^G(S^1)$ fibers over a product of different numbers of circles.   As will become clear, the fibres can be discs, spheres, and more varied spaces. The space $Symm^G_n(S^1)$ can be described explicitly, for all $G$ with $n \leq 4$, and for a wide variety of `useful' groups $G$ with arbitrary $n$. 

\begin{example} \label{ex:orderedTuples} For $G$ the trivial group, $Symm^G_n(S^1) = (S^1)^n$.  This is the space of ordered $n$-tuples of octave-reduced notes. For $n=2$, this is the torus $S^{1} \times S^{1}$.
\end{example}

\begin{example} \label{ex:symm} The space $Symm^{\Sigma_{n-1}}_n(S^1) \cong S^1 \times Symm^{\Sigma_{n-1}}_{n-1}(S^1)$ is a product of $S^1$ and the lower-dimensional symmetric product.  This is the space of ordered pairs, where the first element is an octave-reduced note, and the 2nd element is an unrelated $(n-1)$-note chord.  More generally, $Symm^{\Sigma_i \times \Sigma_j}_n(S^1) \cong Symm^{\Sigma_i}_i(S^1) \times Symm^{\Sigma_j}_j(S^1)$ provided $i+j=n$.
\end{example} 

\begin{example} \label{ex:chordSpace} 
The space of $n$-note chords (under chord-class distance measure (\ref{eqn:distanceCC}))
is the symmetric product of the circle, i.e., $Symm_n(S^1) = (S^1)^n / \Sigma_n$. 
Morton's Theorem \ref{MortonThm} shows that $Symm_n(S^1)$ fibers over $S^1$ with 
fiber a disc, and the bundle is orientable precisely when $n$ is odd.  
There are only two disc bundles over a circle, the trivial (therefore orientable)
one, and the non-orientable one.  The first is a product $S^1 \times D^{n-1}$ while the 
latter can be thought of as a generalized M\"obius band, $S^1 \ltimes D^{n-1}$, and in the case
$n=2$ it is the M\"obius band. 
\end{example}  

\begin{example} \label{ex:s1xs2} 
For $G=A_3 = \Zed_3$, the space $Symm^G_3(S^1)$ represents the 3-note melody space, i.e. the space of triples of octave-reduced notes, where we only remember the cyclic ordering of the notes.  We show in Section \ref{ex:s1xs2} this space is homeomorphic to $S^1 \times S^2$. 
\end{example}

\begin{example} \label{ex:orbifold} The 4-part voice-leading measure of Callender et. al. \cite{callender2005}-\cite{callender2008} is equivalent to the chord-class distance function (\ref{eqn:distanceCC}) and so the symmetric product $Symm_{4}(S^1)$ is the nonorientable orbifold identified in \cite{callender2005}-\cite{callender2008} and explored and popularized by Tymoczko in \cite{Tymoczko2006}-\cite{Tymoczko2010}.
\end{example}

Similar arguments can be used to identify the sub-symmetric products $Symm^G_4(S^1)$.  The subgroup lattice \footnote[1]{We consider the subgroup lattice where the subgroups are taken up to conjugacy.  Congugate subgroups of the symmetric group result in homeomorphic sub-symmetric product spaces.} for $\Sigma_4$ consists of: $A_4$, two different cyclic subgroups of order $2$, $\Sigma_2\times \Sigma_2$ in two different ways (disjoint cycles and not), $\Sigma_3$, $\Zed_4$, and $D_4$.  The corresponding spaces $Symm^G_4(S^1)$ can all be described as disc bundles over products of circles.  
Also, for example, it is possible to explicitly address the topological structure in the case of a $j$-term time-delay embedding of $k$-note chords.
\begin{proposition} \label{prop:delayChords}
When $G = (\Sigma_k)^j$, and $n=kj$, 
\[ Symm^G_n (S^1) \cong (Symm_k (S^1))^j = \prod_j Symm_k (S^1). \]
\begin{proof}
$(\Sigma_k)^j$ can be considered a subgroup of $\Sigma_n$ by partitioning $n =kj$ 
into $k$ sets of size $j$. This is the subgroup of $\Sigma_n$ that preserves the partition. 
There is a function 
\[ Symm^G_n (S^1) \to \prod_j Symm_k (S^1) \]
given by grouping together the $k$-note chords.  It is continuous by the properties of product and quotient space and one-to-one and onto by design.  Since both spaces are compact Hausdorff spaces, this function is a homeomorphism.  

Given two fibre bundles $\pi_i : E_i \to B_i$ for $i=1,2$ the product is naturally a fibre bundle
$\pi_1 \times \pi_2 : E_1 \times E_2 \to B_1 \times B_1$, where $(\pi_1 \times \pi_2)(e_1,e_2) = (\pi_1(e_1),\pi_2(e_2))$. The fibre of the product is $(\pi_1\times \pi_2)^{-1}(b_1,b_2) = \pi_1^{-1}(b_1) \times \pi_2^{-1}(b_2)$ the product of the fibres of the individual maps.  
\end{proof}
\end{proposition}

Since $Symm_k (S^1)$ is a disc bundle over $S^1$ which is orientable only when $k$ is odd, $Symm^G_n(S^1)$ is therefore a disc bundle over $(S^1)^j$ which is orientable if and only if $k$ is odd.  
The case of $Symm^G_4 (S^1)$, where $G=\Zed_4$ is the subgroup of $\Sigma_4$ generated by a $4$-cycle is analogous to Example \ref{ex:s1xs2}, and is similarly involved.  The space $Symm^{\Zed_4}_4 (S^1)$ represents 4-pitch melody space, and Example B.1 presents a brief derivation of its homotopy-type. The space $Symm^{\Zed_4}_4 (S^1)$ can be analyzed via its Morton bundle, as shown in Electronic supplement Appendix B.

\subsection{Finite Subset Spaces} \label{sec:finiteSubspaces}

The $n$-th {\it finite subset space} of a topological space $X$ is the collection of all subsets of $X$ with at least one, and no more than $n$ elements, denoted $exp_n(X)$.   It is given the quotient topology via the map $X^n \to exp_n(X)$ which sends $(x_1,x_2,\cdots,x_n) \longmapsto \{x_1,x_2,\cdots,x_n\}$.  Thus $(0,0,1) = (0,1,1)$ in $exp_3(\{0,1\})$, i.e. the subset space does not keep track of the number of times an element occurs but only keeps track of the binary information that `an element is in the set, or not.'  
This differs from $Symm_3(\{0,1\})$ in that $(0,0,1) \neq (0,1,1)$ in $Symm_3(\{0,1\})$.  
Thus $Symm_n(X)$ can be thought of as being the `multi-sets' with precisely $n$ elements, counting repetition.  
Tuffley studied the finite subset spaces of $S^1$ in his dissertation.  

\begin{theorem}[(Tuffley \cite{Tuff})]  The space $exp_1(S^1)$ can be identified with $S^1$.  
The space $exp_2 (S^1)$ is a M\"obius strip, with $\partial exp_2 (S^1) = exp_1 (S^1)$.
The space $exp_3 (S^1)$ is homeomorphic to $S^3$, moreover if $exp_3 (S^1)$ is visualized as the $3$-sphere, 
$exp_1 (S^1)$ sits in it as a trefoil knot. 
More generally, $exp_n (S^1)$ has the homotopy-type of an odd-dimensional sphere of dimension $n$ or $n-1$ if $n$ is odd or even respectively.  The inclusion $exp_{2n-1} (S^1) \to exp_{2n} (S^1)$ induces multiplication by $2$ on $H_{2n-1}$.  The complement of $exp_{n-2}(S^1)$ in $exp_n (S^1)$ has the homotopy-type of a $(n-1,n)$-torus knot complement. 
\end{theorem}

In particular, note that the quotient map $Symm_3 (S^1) \to exp_3 (S^1)$ is one-to-one 
except on the subspace mapping to $exp_2 (S^1)$.  
Sitting inside $Symm_3 (S^1)$ this is the subspace of the form 
$\{(z_1,z_1,z_2) \}$ where $z_1,z_2 \in S^1$. 
Note that $(z_1,z_1,z_2)$ and $(z_1,z_2,z_2)$ are mapped to the same 
object in $exp_3 (S^1)$.  So $Symm_3 (S^1)$ is equivalent to 
$D^2 \times S^1$ and the map $Symm_3 (S^1) \to exp_3 (S^1)$ is 
one-to-one except on the boundary, where it collapses the torus 
boundary of $Symm_3 (S^1)$ to the M\"obius band $exp_2 (S^1)$. 

\subsection{Homology} \label{sec:homol}

The idea of homology was first introduced by Poincar\'e; a modern treatment can 
be found in \cite{Hatcher}.  
Singular homology is the result of a compromise between the 
geometric intuition of a `hole' and the desire to have something that is readily 
computable.  The $i$-th homology group of a space $X$, $H_i X$, is an abelian group whose 
non-zero elements are objects that detect combinatorial $i$-dimensional `holes' in $X$.  
The precise definition is in the language of simplices and simplicial complexes.  Let $\Delta^n$ be the
standard $n$-simplex, then a {\it singular n-simplex} in a space $X$ 
is a continuous function $\sigma : \Delta^n \to X$.  The {\it boundary} of $\sigma$ is denoted
$\partial \sigma$ and is defined as the formal sum over the bounding $(n-1)$-singular simplices, 
precisely, $\partial \sigma = \sum_{k=0}^n (-1)^k \sigma_{|[0,1,\cdots,k-1,k+1,\cdots,n]}$ where
$[0,1,\cdots,k-1,k+1,\cdots,n]$ indicates we are restricting to the $k$-th face of the $n$-simplex
and parametrizing the face by lexicographically ordering the vertices.  Thus the boundary of
a $1$-simplex is $\sigma_{|1} - \sigma_{|0}$.  In general, this is the orientation assigned
to the boundary as one learns in Stokes' Theorem, in multi-variable calculus.  The {\it group of singular $n$-chains}
in $X$ is denoted $C_n X$ and defined as the group of formal sums of singular $n$-simplices in $X$.  The boundary map
extends to a linear map $\partial : C_n X \to C_{n-1} X$, and the {\em $n$-th homology group} $H_n(X)$ is defined
as the kernel of $\partial :C_n X \to C_{n-1} X$ modulo the image of $\partial : C_{n+1} X \to C_n X$. 
In the case of a simplicial complex, for the purposes of computing homology it is known one can replace
the arbitrary continuous functions $\sigma$ with the characteristic maps of the simplicial complex. Thus
if the simplicial complex is finite, computing the homology of the complex becomes a polynomial-time
computation in the number of input simplices. 
Other work using homology in music theory can be found in \cite{nestke} and \cite{mazzola2012}.


Abelian groups have two standard invariants, the rank and torsion subgroup. 
The {\em torsion subgroup} of an abelian group is the subgroup of elements of finite order, so in 
the group $\Zed^{100} \oplus \Zed_{24}$, the torsion subgroup would be $\Zed_{24}$. 
Taking the abelian group modulo its torsion subgroup (provided this group is finitely-generated), the
minimal number of generators is called the {\it rank} of the original abelian group.  If this group is
not finitely generated, the group is said to have infinite rank.  
Hence the group $\Zed^{100} \oplus \Zed_{24}$ has rank $100$. 

The rank of the $i$-th homology group of $X$, $H_i X$, is called the {\em $i$-th Betti number} of $X$, 
and is typically denoted $\beta_i(X)$.   For example, the only invariant of $H_0 X$ is $\beta_0 (X)$, and
this is the number of path-components of $X$. The number of
linearly independent loops in $X$ that do not bound discs is $\beta_1 (X)$.  When the space is 
a circle, $H_1 S^1 \cong \Zed$ so there is one independent loop.  
Similarly, $H_1 (S^1 \times S^1) \cong
\Zed^2$ and there are two linearly independent circles in the torus, correspondingly
$\beta_1(S^1\times S^1) = 2$.  Homology can be somewhat subtle. If
$X$ is the Klein bottle, $H_1 X \cong \Zed \oplus \Zed_2$, there is a
torsion subgroup corresponding to the circle in the Klein bottle that does
not reverse orientation; in this case $\beta_1(X) = 1$ since Betti numbers 
ignore torsion.

\subsection{Persistent Homology} \label{sec:persHomology}

Introduced by Carlsson and his coworkers \cite{Carlsson} as a way of parsing large data sets, persistent homology 
takes as its starting point the functoriality of homology.  Precisely, if $m : X \to Y$ is a continuous function, there is an induced homomorphism of homology groups $m_* : H_i X \to H_i Y$, given by composing the singular simplices in $X$ with $m$.  This allows a comparison between the homology groups of two spaces. 

Let $P$ be a collection of points in a space $X$ where there is a notion of 
distance $d : X \times X \to \mathbb R$. The {\it singular homology groups} $H_* P_\epsilon$
are used to represent the general shape of the $\epsilon$-neighborhood of $P$ in $X$, 
where $P_\epsilon = \{ x \in X : d(x,p) < \epsilon \text{ for some } p \in P\} \subset X$.  


Figure \ref{fig:pointcloudfig} shows a set $P$ consisting of 21 points and the
$\epsilon$-neighborhoods for five cases. For $\epsilon=0.6$, the only non-trivial homology group is
$H_0 P_{0.6} \cong \Zed^{21}$.  The top-middle figure
shows $\epsilon=0.9$, and $H_0 P_{0.9} \cong \Zed^{12}$ is the only non-trivial homology group.
In the third figure, all homology groups of $H_0 P_{1.6}$ are trivial 
(the homology groups of a point). As $\epsilon$ increases yet further to 2.4, a single
large circle appears, generating a $1$-dimensional homology class. For still larger $\epsilon$ (such as 3.5), $P_\epsilon$ is approximately a disc, and all homology groups disappear, with the exception of a single $H_0$ class. 

\addtocounter{figure}{2}
\begin{figure}
\begin{center}
\includegraphics[width=3.5in]{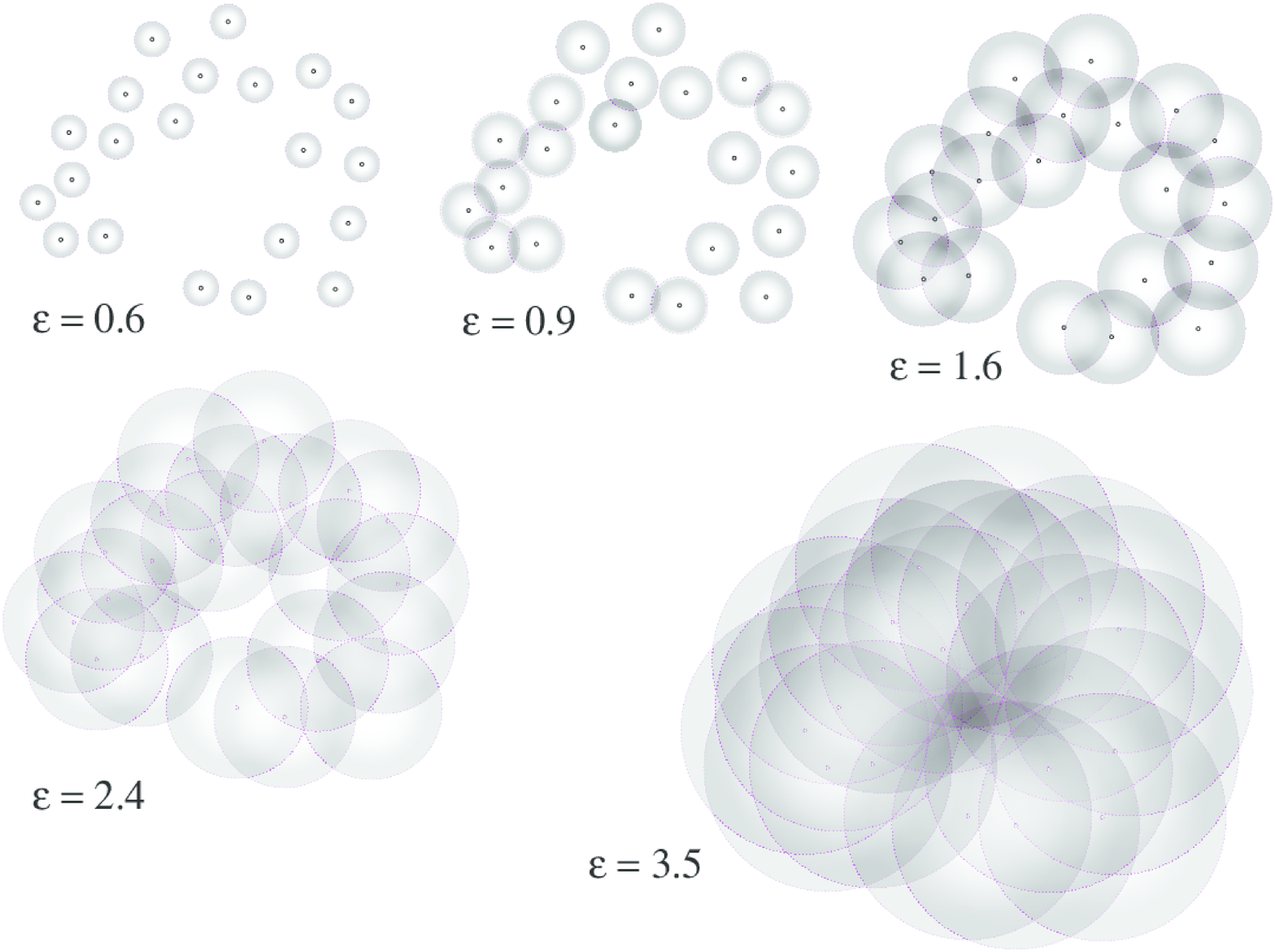}
\end{center}
\caption{Representation of a point cloud of $21$ points and $P_\epsilon$ for $\epsilon=0.6, 0.9, 1.6, 2.4$, and $3.5$}
\label{fig:pointcloudfig}
\end{figure}

The {\em barcode} of a point cloud $P$ in $X$ is a record of how the homology
of the epsilon-neighborhood $P_\epsilon$ varies as a function of $\epsilon$.  
If $\epsilon_1 < \epsilon_2$, $P_{\epsilon_1} \subset P_{\epsilon_2}$
so there is an induced map of homology groups $H_i P_{\epsilon_1} \to H_i P_{\epsilon_2}$. 
If a homology class in $H_i P_{\epsilon_1}$ is mapped to a non-trivial homology
class in $H_i P_{\epsilon_2}$, it is said to {\it persist} over the interval $[\epsilon_1,\epsilon_2]$
and the {\it bar} associated to a homology class is the largest interval over which it
persists.  The {\it barcode} of the point cloud $P$ is the record of all the bars
over a range of $\epsilon$.  Reading a barcode from left-to-right over the 
$\epsilon$-axis, bars may appear and disappear. 
A new bar means that a new homology class appears at that $\epsilon$; the bar remains for as long as that 
homology class persists, and then it vanishes. For example, in Figure \ref{fig:pointcloudfig},
consider the seven left-most points. For small $\epsilon$, these are isolated and 
each appear as a separate line in the barcode. By $\epsilon=0.9$, these seven lines
merge together into a single line. All 21 lines merge into one by $\epsilon=1.6$ and
a single $H_{0}$ line persists for all larger $\epsilon$.
Figures \ref{fig:kingsStandardRhythm1D}, \ref{fig:majorScale}, \ref{fig:abbotts1D},
\ref{fig:abbotts2D} etc. provide many examples of barcodes.

Unfortunately, as described above, the computation of persistent homology requires 
knowledge of more than just the point data $P$.  The intermediate step of considering
$P_\epsilon$ requires knowledge of the distance function on the entire space $X$, which 
can be difficult to compute. A simplicial complex 
$P^\Delta_\epsilon$ (called the {\it Vietoris-Rips complex}),
which depends only on the distance relation between individual points in $P$, 
can be used to approximate $P_\epsilon$. The simplicial complex $P^\Delta_\epsilon$ is constructed inductively. $P$ is the vertex-set of $P^\Delta_\epsilon$.  For every pair of points of $P$ separated by less than $\epsilon$ one adds an edge to the simplicial complex.  Similarly, for every $(n+1)$-tuple of points of $P$ which are pairwise distance less than $\epsilon$ apart, one adds an $n$-simplex to the simplicial complex. The simplicial complex $P^\Delta_\epsilon$ is the result of this inductive procedure. Thus persistent homology should be viewed as a coarse measure of the shape of the data $P$ {\em as an abstract space}.  
Theorem 2.3 of Adler et. al. \cite{larry} gives
a sufficient criterion for $H_* P_\epsilon$ and $H_* P^\Delta_\epsilon$ to be
isomorphic, although the preconditions of their theorem almost never hold 
with the data considered in this paper. 

For $\epsilon$ small, the structure is always the same; $P^\Delta_\epsilon = P$, as all 
points are separated from all other points. For $\epsilon$ larger than the diameter of $P$, 
the structure is simple to describe; $P^\Delta_\epsilon$ consists of a single simplex whose 
vertex set is $P$.  In between small and large $\epsilon$ there may be a wide
variety of homology groups.  Homology classes that persist over a significant $\epsilon$-interval
are called {\em persistent} and frequently reflect interesting underlying structure
in the data set $P$.  A number of applications have begun to appear in areas such as 
image processing \cite{Carlsson} and in the analysis of biological data \cite{singh}. 

\section{Simulations with Musical Data} \label{sec:simulations}

This section analyzes a number of musical pieces (some contrived and some from 
various standard repertoires) using persistent homology
by turning raw musical data (in this 
case, note-level MIDI data) into a collection of points. 
Using the appropriate 
distance measures (\ref{eqn:distanceRhy})-(\ref{eqn:distanceCyc}), 
it is straightforward to calculate the distance between
every pair of points in the data. 
The {\tt javaPlex} software \cite{javaPlex}, designed to 
``calculate the persistent homology of finite 
simplicial complexes... generated from point cloud data,''
can be used to calculate the barcodes for a given set of data.
See Section \ref{sec:persHomology} for the basics of persistent homology, and
how to read the barcodes.
The input is a set of distances (between all
points in the point-cloud) and the output is a set of plots that show
the Betti numbers as a function of the persistence
parameter $\epsilon$. 

\subsection{Rhythmic Barcodes} \label{sec:rhythm}

The first example translates the ``Ewe'' rhythm of 
Figure \ref{fig:yoruba} into the set of time points
\begin{equation} \label{eqn:ewerhy}
\left\{ 0, \frac{2}{12}, \frac{4}{12}, \frac{5}{12}, \frac{7}{12}, \frac{9}{12},
\frac{11}{12} \right\} .
\end{equation}
The distance is calculated between every pair of these time points 
under the necklace metric (\ref{eqn:distanceRhy}), and this 
set of distances is the ``point cloud'' that is input into the javaPlex software.
The resulting barcodes are shown in Figure \ref{fig:kingsStandardRhythm1D}.

\addtocounter{figure}{3}
\begin{figure}
\centering \includegraphics[width=4.8in]{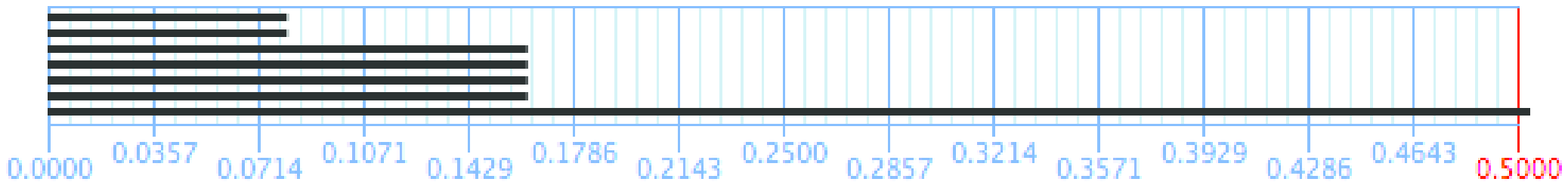}
\includegraphics[width=4.8in]{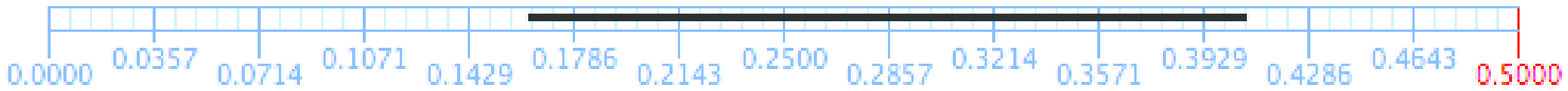}
\caption{Barcodes calculated by the javaPlex software show the
number of connected components in the dimension 0 plot (top) and
the number of circles in the dimension 1 plot (bottom), as the size
parameter $\epsilon$ varies. 
Barcodes for King's ``standard rhythm'' \cite{king61} of Figure 
\ref{fig:yoruba} show the distribution of time intervals in the
rhythm in the top (dimension 0) plot and show the circular structure
with Betti numbers $\beta_{0}=1$ and $\beta_{1}=1$ in the bottom (dimension 1) plot.
Using the distance function (\ref{eqn:distanceRhy}),
the largest distance between elements is 0.5, and the diameter of the
space, homeomorphic to $S^{1}$, is also 0.5.
\label{fig:kingsStandardRhythm1D}}
\end{figure}


The dimension 0 barcode shows the clustering of the points in time.
In the sequence (\ref{eqn:ewerhy}), the minimum distance is 
$\frac{1}{12}$, and this occurs in two places, between the 
third and fourth elements of the set, and again between the 11th and the first elements.
Accordingly, the barcode shows two lines that vanish when $\epsilon$
reaches $0.08$. Since the largest distance between any two adjacent
time points is $0.16$, all the points merge into one cluster
at $\epsilon=0.16$. To see why this happens, imagine that a balloon
is placed around each of the small circles in Figure \ref{fig:yoruba}.
When the radius of the balloon is small, each balloon is isolated and contains 
no other circles than the one at its center.
But as the balloons are inflated (which corresponds to increasing $\epsilon$),
the balloons around point $a$ and $f$ grow large enough to touch the points
$b$ and $e$. This is the radius at which the first two lines in the
upper figure disappear, the points $a$ and $b$ are now encompassed 
by the same balloon (and similarly for the points $e$ and $f$). 
As the balloons inflate still further, they grow to include other points that are
two time units apart; at this size, all adjacent points are now connected.
This situation is shown in the dimension 1 barcode which displays a persistent 
$\beta_{1}$ bar from $0.16<\epsilon<0.408$. This is the anticipated
cycle around the necklace. With $\beta_{0}=1$ and $\beta_{1}=1$,
the complete structure is indeed the circle $S^{1}$.
While this is not surprising, it is useful to verify that the 
simulation method and the analysis agree.

All three of the rhythms indicated in Figure \ref{fig:yoruba}
have identical barcodes because the rhythms are circular shifts of each other
and hence contain the same collection of rhythmic intervals.

\subsection{The $C$-Major Scale} \label{sec:notes}

\begin{figure}
\centering \includegraphics[width=4.8in]{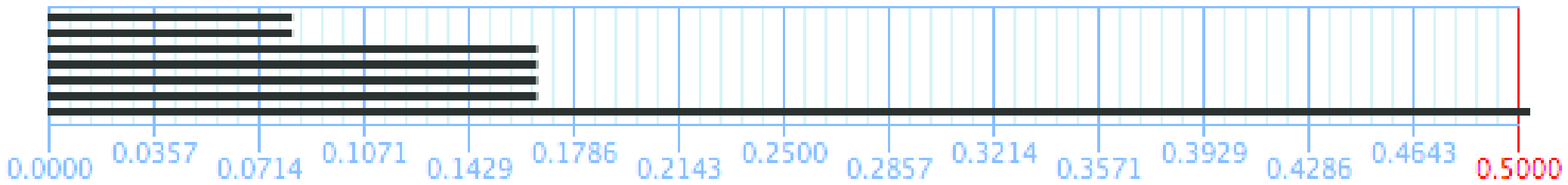}
\includegraphics[width=4.8in]{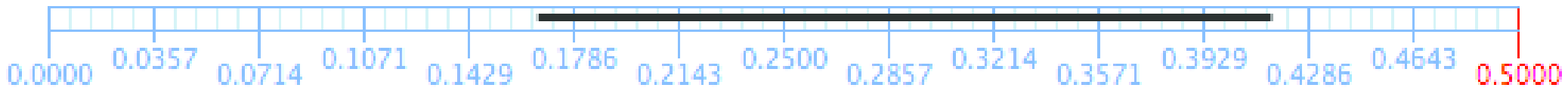}
\caption{Barcodes for the $C$-major scale show seven distinct
lines in the top (one corresponding to each note) for small
$\epsilon$. Two half steps merge at $\epsilon \approx 0.08$
and the whole steps merge at $\epsilon \approx 0.16$, at which
point the circle $S^{1}$ with $\beta_{0} = \beta_{1} = 1$
appears. The distance function used is (\ref{eqn:distancePC})
and the diameter of the space $S^{1}$ is again 0.5. 
\label{fig:majorScale}}
\end{figure}

To verify that the distance function (\ref{eqn:distancePC}) makes sense, 
consider the $C$-major scale 
consisting of the eight notes $C, D, E, F, G, A, B, C$ with 
frequencies given starting at $C=261.6$ Hz.
The data is defined by a matrix of 
distances between all pairs of the eight notes using the 
pitch-class distance on $\mathbb{R}_+$.
The resulting barcodes are shown in Figure \ref{fig:majorScale}. 

These two plots are straightforward to interpret. When the size parameter 
$\epsilon$ is small, there are seven distinct notes. Though we
input eight notes, the high $C$ has exactly the same distances 
to all the other notes as the low $C$ under the pitch-class distance
(\ref{eqn:distancePC}). Since the distance from the high $C$ to the 
low $C$ is zero, the barcode
merges these two tones even at $\epsilon=0$. When $\epsilon$ reaches $0.08$, 
the two half steps (the intervals between $E$-$F$ and 
$B$-$C$) merge. When $\epsilon$ reaches $0.16$, the five
remaining connected components (all the major seconds)
merge into one. Thus $\beta_{0}= 1$ for all greater $\epsilon$.
At $\epsilon = 0.16$, the dimension 1 code shows a single 
component, which persists until $\epsilon = 0.4$. This $\beta_{1}= 1$ 
feature is $S^{1}$, and corresponds to the circular structure illustrated 
in Figure \ref{fig:circle-of-notes}.

\subsection{Pitch-Class Barcodes: Abbotts}

\begin{figure}
\centering \includegraphics[width=4.8in]{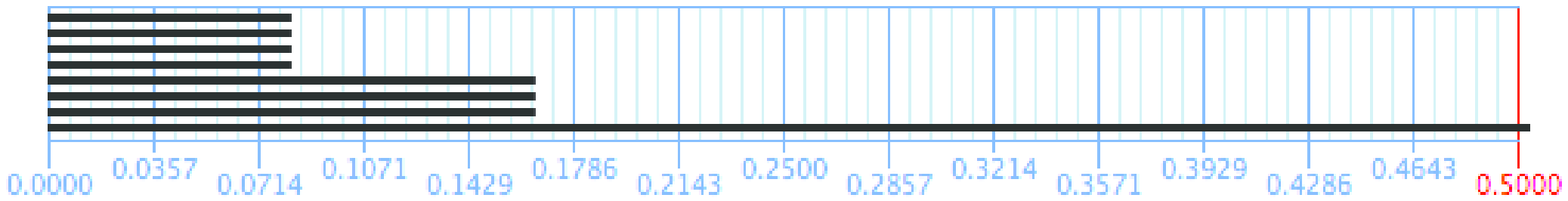}
\includegraphics[width=4.8in]{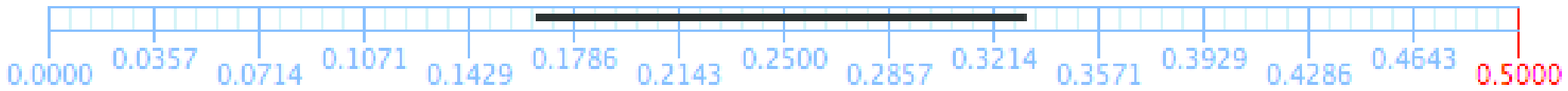}
\caption{Barcodes for the traditional folk tune ``Abbott's Bromley
Horn Dance'' (see Figure C1) 
under the distance function (\ref{eqn:distancePC}) show many of the
same features as the major scale barcodes of Figure \ref{fig:majorScale}.
The distribution of whole and half steps are clear from the 
$\beta_{0}$ code for small $\epsilon$ while
$S^{1}$ appears again in $\beta_{1}$ for $0.16<\epsilon<0.33$.
\label{fig:abbotts1D}}
\end{figure}

The input data for example of Figure \ref{fig:majorScale} was built specifically
with the circular structure of Figure \ref{fig:circle-of-notes} in mind, 
so it is perhaps unsurprising that the circle appears. Will such shapes appear in real music?
The website \cite{peterson} contains a large selection of 
traditional melodies, with most tunes available in both sheet music
and as standard MIDI files. The musical score for ``Abbott's Bromley
Horn Dance'' is shown in Figure C1 and
the corresponding barcodes are shown in Figure \ref{fig:abbotts1D}.
The distance between note pairs is given by (\ref{eqn:distancePC}).

The top barcode in Figure \ref{fig:abbotts1D} shows eight lines,
which correspond to the eight notes that appear in the score
(observe again the insensitivity to octave). Four disappear at 
$\epsilon = 0.08$, which correspond to the four half steps 
($F\sharp$-$F$, $D\sharp$-$E$, $B$-$C$, and $D$-$D\sharp$). 
Three more disappear at $\epsilon = 0.16$. Along with the 
constant bar, these correspond to the four whole steps ($E$-$F\sharp$, $G$-$A$,
$A$-$B$, and $C$-$D$). All of these join
into one bar for all larger $\epsilon$.
The region $0.16<\epsilon<0.33$ is characterized by 
$\beta_{0}=1$ (one connected component) and $\beta_{1}=1$,
one circle. This is again the circular structure $S^{1}$.
In fact, all the melodies from the website \cite{peterson}
show this same structure, though the number of half and whole
steps changes to reflect the scale of the piece, and the exact
extent of the $\beta_{0}= \beta_{1}=1$ region is somewhat variable.

\subsection{Barcodes for an Ussak Makam}

\begin{figure}
\centering \includegraphics[width=4.8in]{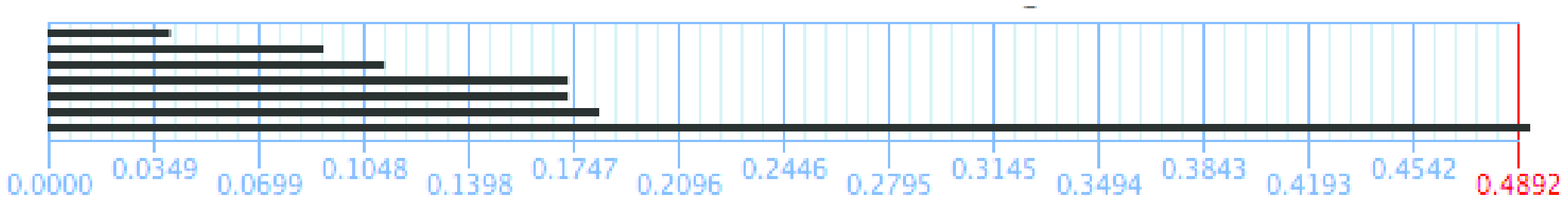}
\includegraphics[width=4.8in]{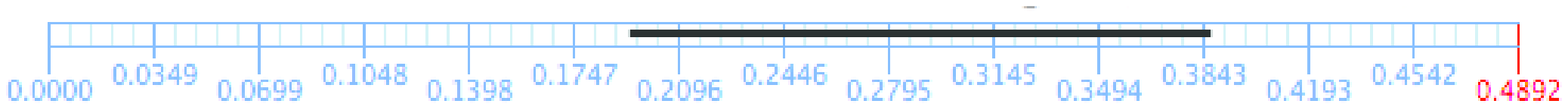}
\caption{Barcodes for an Ussak makam  
under the distance function (\ref{eqn:distancePC}) show many of the
same features as the barcodes of Figure \ref{fig:abbotts1D}, though
some significant dissimilarities are also apparent.
The occurrence of many different intervals is apparent in the 
staggered appearance of the dimension 0 barcode (top).
Again, the circle appears in the bottom barcode for 
$0.18<\epsilon<0.38$.
\label{fig:ussak1D}}
\end{figure}

While it might appear from the above examples that the analysis requires
the pitches to be subsets of the semitones of the Western equal-tempered scale,
they can assume any possible values. The Turkish makam tradition utilizes
many tones/pitches that are distinct from the Western pitches, and 
it is instructive to look at the bar codes for one such piece.
A data base of 1700 makams available online \cite{kemal} contains MIDI and 
score data. About a hundred of the makams are from the mode ussak, and these
were analyzed using the javaPlex software. The barcodes for a typical example
are presented in Figure \ref{fig:ussak1D}. Like the Western pieces, 
as $\epsilon$ increases, the pitches merge together. Unlike the equal tempered 
pieces, however, there are many different intervals between various pitches 
and so the bars of the code merge together at somewhat staggered values
(rather than at one or two discrete $\epsilon$ values). 
Nonetheless, the structure of the underlying space is revealed 
by the merging of all the pitches in the dimension zero barcode by 
$\epsilon \approx 0.18$ and the emergence of the dimension 1 bar at
about this same value. Again, this represents the circular $S^{1}$ structure.

\subsection{Time-Delay Embedding: 2-D Abbotts} \label{timeDelay1}

\begin{figure}
\centering \includegraphics[width=4.8in]{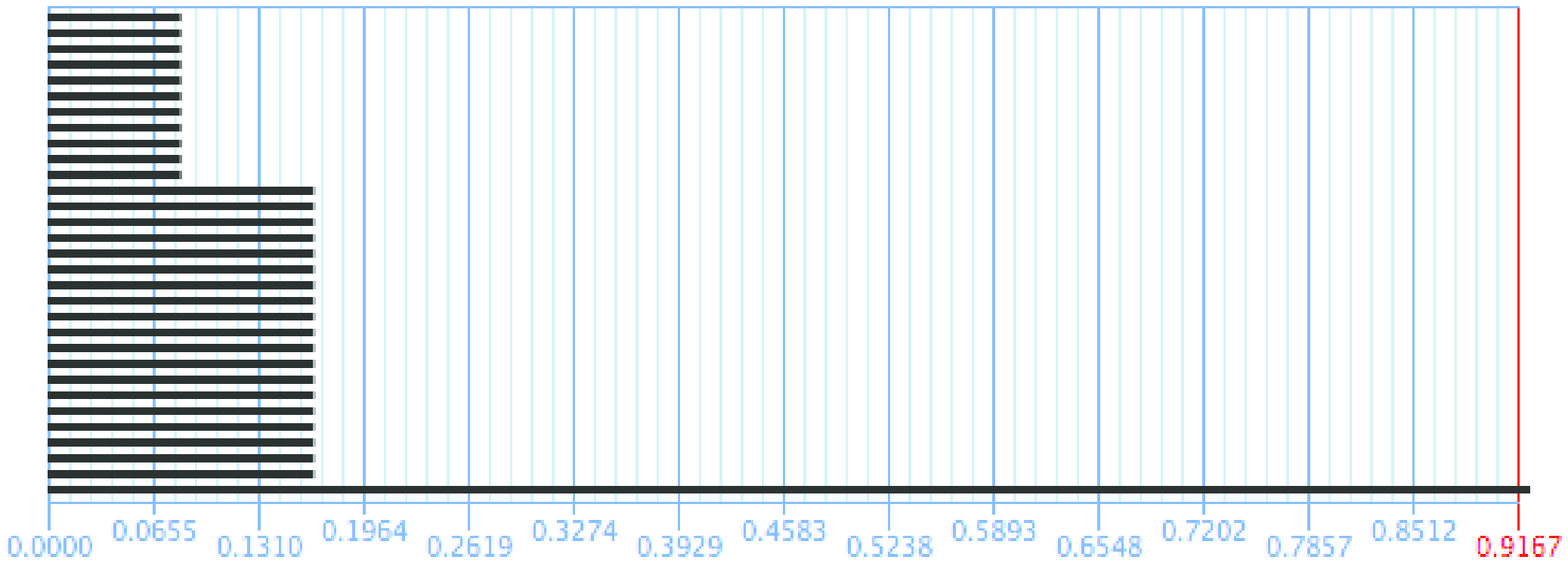}
\includegraphics[width=4.8in]{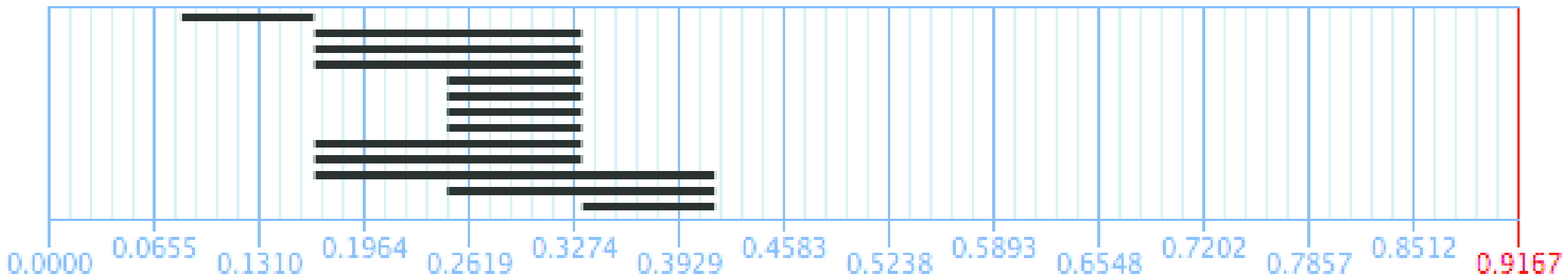}
\includegraphics[width=4.8in]{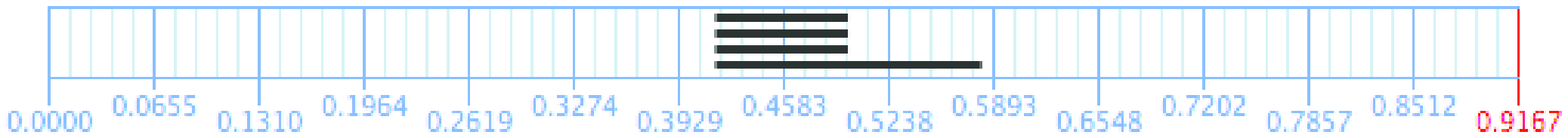}
\includegraphics[width=4.8in]{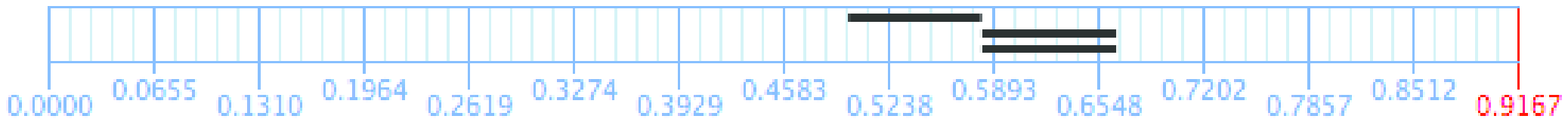}
\includegraphics[width=4.8in]{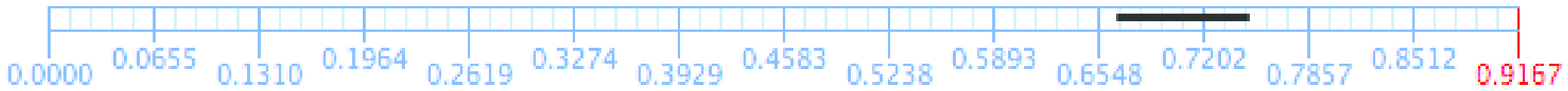}
\caption{Barcodes for the two-dimensional time-delay 
embedding of ``Abbott's Bromley Horn Dance'' (see Figure C1) 
are considerably more interesting than those in Figure \ref{fig:abbotts1D}
without the time-delay.
The top picture shows the dimension 0 barcodes, the second shows the 
dimension-one barcodes, continuing through the dimension-four barcodes.
The distance function (\ref{eqn:distancePC}) is applied to 
each element of 2-vectors (which consist of consecutive pairs of note values), 
and then summed. The topological space in
which this is embedded is $S^{1}\times S^{1}$, the 
$n=2$ case of Example \ref{ex:orderedTuples}.
The maximum distance between elements of the performance is 
$0.9167$, which is close to the diameter of the complete space.
\label{fig:abbotts2D}}
\end{figure}

The analyses of Figures \ref{fig:majorScale} and \ref{fig:abbotts1D}
may be somewhat naive because they suppress temporal 
information in the melody. A time-delay
embedding considers successive pairs (or more) of elements and so 
recovers some of the temporal relationships; it is a common procedure in time series analysis.
For this approach the appropriate distance function is
(\ref{eqn:distancePC}). 
Building a matrix of all such distances for ``Abbott's Bromley Horn Dance''
and calculating the barcodes gives Figure \ref{fig:abbotts2D}.

The dimension $0$ barcodes (the top plot in Figure \ref{fig:abbotts2D})
can be interpreted as showing the distances between consecutive pairs of notes
as the melody progresses over time. Thus there are 11 pairs of notes
that are at a distance of one-half step, since 11 lines end at
$\epsilon=0.08$. There are 19 pairs that differ by a whole step
since 19 lines end at $\epsilon=0.16$. Above this value, all pairs
have merged into a single connected component. This can be interpreted
as saying that the melody progresses primarily by stepwise motion,
and that no pairs of tones are isolated from any other pairs of tones
(though of course there are many individual pairs with larger distances).

The dimension $1$ barcodes (the second plot in Figure \ref{fig:abbotts2D})
shows the number of independent circles present at each value of $\epsilon$, 
the dimension $2$ barcodes (the third plot) show the number of 
closed oriented 2-cycles (such as hollow spheres) 
as a function of $\epsilon$, and the dimension 
$3$ barcodes (the bottom plot) show the distribution of 4-cycles.  
For this example, the embedding space is the torus $S^1 \times S^1$, as in
Example \ref{ex:orderedTuples}. 

This is our first encounter with barcodes beyond the dimension of the underlying
metric space, which is two for the space $S^1 \times S^1$. In brief what we see is a decreasing
number of bars (as the dimension increases) and the persistence intervals of these
bars appears to shrink as the dimension increases.  This is a typically-observed 
phenomenon with data that respects a simple distribution, like a Gaussian \cite{balzano}.  More precisely, what
one typically sees in the barcodes of a Gaussian field is a void supporting no homology in the centre of the distribution, while homology classes get sparser (and higher dimensional) as one approaches the void from the outside.  In this case, there is one 1-dimensional homology class that persist over a fairly wide interval, indicating the data diverge somewhat from this expectation.  In this case it is due to the relative diversity of two-note sequences starting with an $E$, 
and the relative sparsity of other combinations.  

The higher dimensional barcodes are few and have relatively short persistence intervals so they 
represent relatively insignificant features of the data. Moreover, the fact that they represent
cycles that can not exist in the parameter-space $S^1 \times S^1$ further indicates the
idiosyncratic nature of these cycles. This is reinforced by analyses of other songs
from the same set (not shown here) of 2D time-delay embeddings which have analogous 
(but different) groupings of higher order barcodes.

\subsection{Time-Delay Embedding: 3-D Abbotts} \label{timeDelay2}

Longer temporal information can be incorporated by using longer time-delay
embeddings. If a melody consists of a sequence of notes with pitches 
at $f_{1}, f_{2}, f_{3}, f_{4}... $, these can be 
combined into three-tuples (a three-dimensional time-delay
embedding) by creating the sequence $(f_{0}, f_{1}, f_{2})$, 
$(f_{1}, f_{2}, f_{3})$, $(f_{2}, f_{3}, f_{4}), ...$.
The distances between such three-tuples can be calculated by 
adding the distances between the notes element-wise 
using the pitch-class distance. 
Building a matrix of all such distances for ``Abbott's Bromley Horn Dance''
and calculating the barcodes gives Figure C2.

Again, it is straightforward to interpret the dimension $0$ 
barcodes as distances between three-tuples of notes in the melody.
Since the embedding here is $S^1 \times S^1 \times S^1$ by Example \ref{ex:orderedTuples},
the higher dimensional structures are again linked to the pieces
being analyzed.

\subsection{Barcodes for the Circle of Fifths} \label{sec:chords}

\begin{figure}
\centering \includegraphics[width=3.0in]{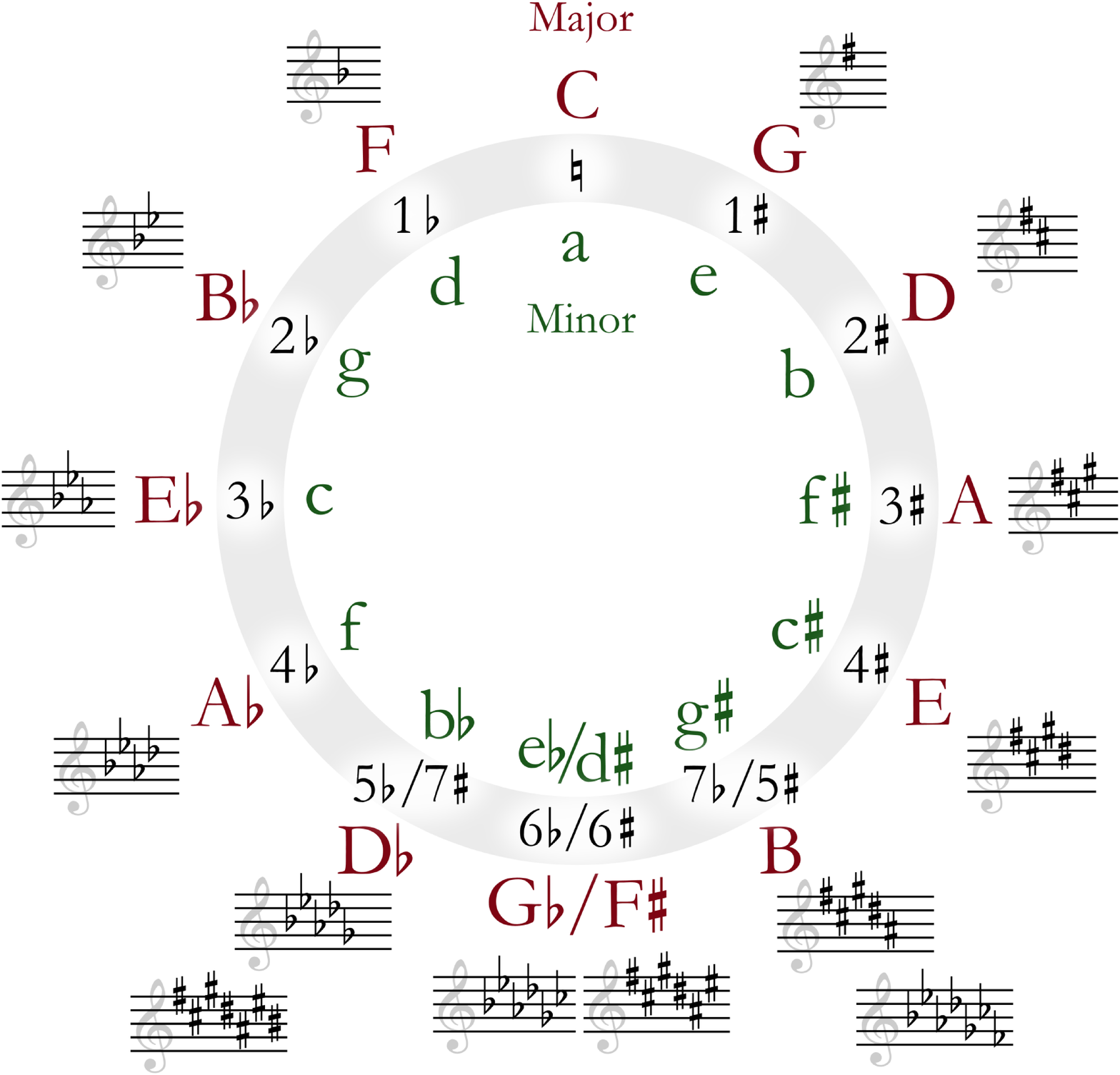}
\caption{The circle of fifths shows the relationships among the 
tones of the equal tempered chromatic scale, standard key signatures, 
and the major and minor keys.  
\label{fig:circle-of-fifths}}
\end{figure}

A well known structure in music theory is the
``circle of fifths,'' shown here in Figure \ref{fig:circle-of-fifths},
which is taken from the Wikipedia article of the same name \cite{wikicircleoffifths}.
The circle of fifths is a standard way musicians and composers talk 
about the close relationships between musical scales and keys,
and represents another way of interpreting the distance between
musical chords and scales. 
The appropriate distance function is (\ref{eqn:distanceCC}), and 
it is reasonable to ask if the circle of fifths can be
recovered from a data analysis. 

\begin{figure}
\centering \includegraphics[width=4.8in]{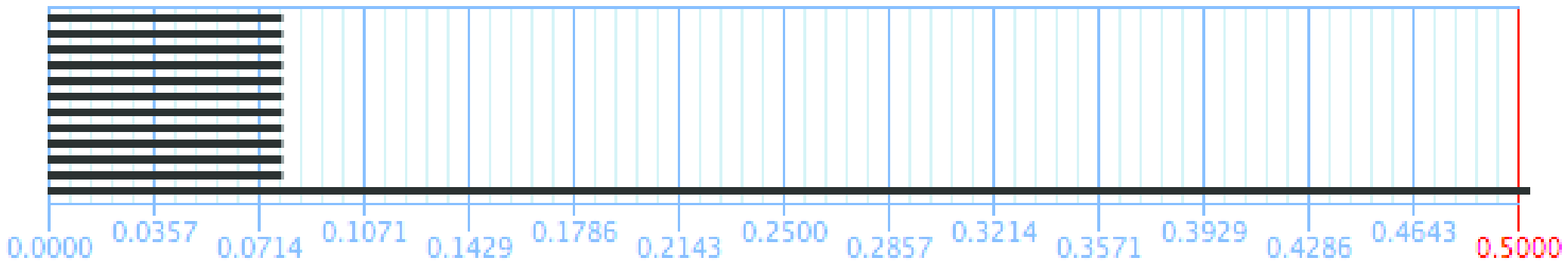}
\includegraphics[width=4.8in]{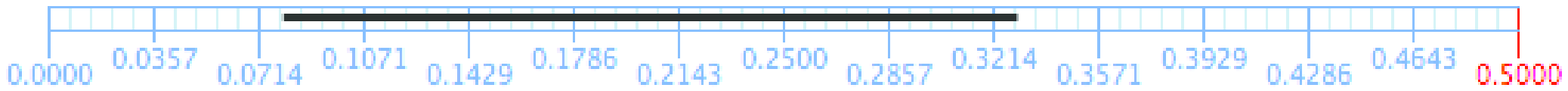}
\includegraphics[width=4.8in]{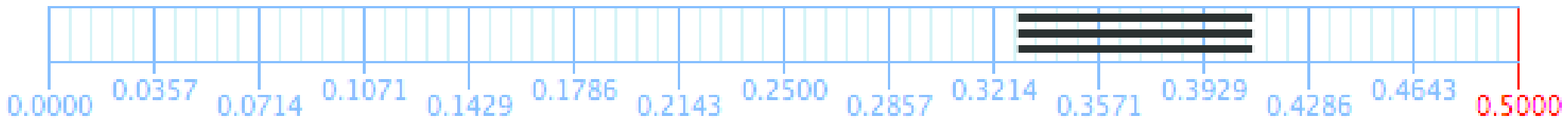}
\caption{Barcodes for a sequence consisting of
a major scale proceeding around the circle of fifths. The ``circle'' 
is the persistent line from $0.08<\epsilon<0.33$
in the dimension 1 plot. Under the chord-class distance function
(\ref{eqn:distanceCC}) the space is $S^{1} \times D^{2}$ 
(this is the $n=3$ case of Example \ref{ex:chordSpace})
and the total diameter of the space is 3.5 while  
the maximum distance between any two elements is 0.5.
\label{fig:circleofFifths}}
\end{figure}

Consider a progression that moves around the circle of fifths:
$C$ major to $G$ major to $D$ major etc, all the way
back to $F$ and finally $C$. For the present example, 
each key is represented by its major scale,
a seven note set of pitches.
Inputting these sets into javaPlex and calculating the barcodes gives Figure
\ref{fig:circleofFifths}. Under the distance (\ref{eqn:distanceCC}), scales that are
a fifth apart (such as $C$ major and $G$ major) have 
a distance of $0.08$. This explains the twelve lines
that merge down to a single connected set at $\epsilon=0.08$
in the dimension $0$ (top) plot.
For $0.08<\epsilon<0.33$, the dimension 1 barcode
shows a single persistent bar; this is precisely the circle of
fifths.

As shown in Example \ref{ex:chordSpace},
the space in which scale (and/or chord) progressions lie 
is an $n$-fold symmetric product of the circle.
As constructed in this section, the circle of fifths can then be 
interpreted as twelve points in $Symm_7(S^1)$.
Observe that the circular structure of Figure \ref{fig:circle-of-fifths}
is not used in any way to draw the barcodes: rather, a set of scales
(along with a metric on those scales) are the only inputs to the software. 
The barcodes then exhibit a structure that can be interpreted as identical 
to the circle.

There are also some higher dimensional features for larger $\epsilon$.   The 
circle of fifths has dihedral symmetry -- specifically, the symmetry group of 
a regular $12$-sided planar polygon (a {\em dodecagon}).  Interestingly the 
metric is not that of a regular dodecagon.  In this case, the distance function is like
the seam of a baseball -- stretching off course in places, then turning back 
later.  Precisely, the $n$-th and $(n+4)$-th note are within $0.33$ of each other {\it for all $n$}.  So the triple consisting of the $n$-th, $(n+4)$-th and $(n+8)$-th notes form a triangle in the Vietoris-Rips complex for each $n$, provided $\epsilon > 0.33$, and they {\it do not} provided $\epsilon < 0.33$.  Since there are $12$ notes, that gives $4$ triangles being attached to the Vietoris-Rips complex near $\epsilon = 0.33$. Attachment of any of these triangles
removes the only $H_1$-class in the Vietoris-Rips complex for $\epsilon < 0.33$.  This is because the
generator of the $H_1$-class is the sum of the edges consisting of consecutive chords in the 
cyclic order of the circle of fifths.  Moreover, the edge consisting of the $n$-th note and the
$(n+4)$-th note is readily computed to be a boundary.  Thus, all four of the triangles remove the
same $H_1$-class from the chain complexes for $\epsilon < 0.33$, so the union of the four triangles
not only removes the $H_1$-class but contributes three $H_2$-classes, via a direct Meyer-Vietoris
argument. 

\subsection{Bar Codes for Chords in Bach's Chorale No. 19}

The musical score to Bach's Chorale No. 19 is
shown in Figure D1.
A MIDI file of this piece, from \cite{classicalMidiArchives},
is parsed to extract the four voices. The distances between all
four-part chords are calculated according to the chord-class
distance (\ref{eqn:distanceCC}), and the results are input to the javaPlex software
in order to draw the barcodes, which are shown in Figure D2.
The dimension 0 barcode shows a large number of chords that are
separated by $\epsilon=0.08$, a somewhat smaller number of chords that 
are separated by a distance of $\epsilon=0.16$, and two 
chords separated by $\epsilon=0.23$. Above this value, 
all chords merge into one connected component.

The dimension 1 barcode in Figure D2
shows $\beta_{0}=6$ connected components
and one circle $\beta_{1}=3$ for $0.16<\epsilon=0.24$, and this 
structure then changes to $\beta_{0}=1$ and $\beta_{1}=3$
for $0.24<\epsilon=0.33$. Features such as these appear to be unique
identifiers of the particular pieces, meaning that other
Bach Chorales from the same data set have different 
Betti numbers that occur over different ranges of $\epsilon$.

\subsection{Time-Delay Embedding of Chords in Abbotts}

It is also easy to include temporal information with a time-delay 
embedding of the chords. Some care must now be taken with the 
distance function since it is a mixture of the chord-class
(\ref{eqn:distanceCC}) and the sum-of-elements distance useful 
in time-delay embeddings as in Section \ref{timeDelay1}
and \ref{timeDelay2}. 
To be explicit, suppose that a sequence
of 3-note chords
$F=(f_{1}, f_{2}, f_{3})$, $G=(g_{1}, g_{2}, g_{3})$, $H=(h_{1}, h_{2}, h_{3})$,
$I=(i_{1}, i_{2}, i_{3}), ...$
is gathered into time-delay pairs to form 6-tuples
$(F, G)$, $(G, H)$, $(H, I)...$. 
Then the distance between pairs of such
6-tuples is the time-delay embedding distance
\begin{equation} \label{eqn:distanceTDE}
d_{TDE}((F,G),(H,I)) = d_{CC}(F,H)+d_{CC}(G,I).
\end{equation}
Observe that under this measure, the distance between notes within a single chord
is measured differently from the distance between the same notes when they occur
at different times.
An example of the resulting barcodes is shown in 
Figure C1 for a 2D time-delay embedding of 3-note chords. 
The barcode shows between several circles for 
$\epsilon$ in the region of 0.5 indicating some structure in the 
chord pattern. The underlying space is described in Proposition \ref{prop:delayChords} 
for the case of a $j=2$-term time-delay embedding of $k=3$-note chords.

\subsection{Two-Dimensional Time-Delay Embedding of Chorale No. 19}

Figure D3 shows the two-dimensional
time-delay embedding of the Bach's Chorale No. 19, where the zero
dimensional plot is interpreted directly in terms of the distribution of
chord pairs and how they cluster under the chord-class distance.
Again, there is a collection of persistent $\beta_{1}$ 
circles and scattered entries in the dimension 2 and 3 barcodes.

\subsection{Rhythms and Cyclic Permutations}

Repetitive rhythmic patterns with period $n$
may be written in the timeline notation as elements of $\Zed^{N}$ 
where each ``1'' represents the occurrence of a pulse and 
each ``0'' represents silence. 
For instance, the 16-beat {\it clave son} can be transcribed as
\begin{equation} \label{eq:claveson}
\small{\mathbf{Clave \ Son}} = \left( 1,0,0,1,0,0,1,0,0,0,1,0,1,0,0,0 \right) \in \Zed^{16}.
\end{equation}
Table \ref{tab:afro-cuban} shows a variety of Afro-Cuban rhythms
each of which is the basis of a traditional or popular musical style.
From a musical perspective, it is clear that some such rhythms can be 
thought of as ``more alike'' than others:
they may be performed in similar musical styles, they may 
occur close together geographically, or there may be
cultural ties between groups that prefer certain patterns
over others. It is intuitively plausible that rhythms may be ``close together''
or ``far apart'' based on similarities and differences in the set of 
temporal intervals that define the rhythms. 
An appropriate notion of distance in this case is 
(\ref{eqn:distanceCyc}), which considers all patterns related by a cyclic shift 
to be identical.  In this case, the underlying ambient space would be
$Symm^{\Zed_5}_5(S^1)$, a space consisting of $5$-tuples 
that represent onset times in a temporal cycle.  Technically, the rhythm
space should be the subspace of distinct $5$-tuples, but $Symm^{\Zed_5}_5(S^1)$ has
the correct metric on this subspace so we create no special notation for this situation.

\begin{table}
\caption{Some Afro-Cuban rhythms}
\begin{center}
\begin{tabular}{lcccccccccccccccc} 
bossa-nova & 1 & 0 & 0 & 1 & 0 & 0 & 1 & 0 & 0 & 0 & 1 & 0 & 0 & 1 & 0 & 0 \\
gahu       & 1 & 0 & 0 & 1 & 0 & 0 & 1 & 0 & 0 & 0 & 1 & 0 & 0 & 0 & 1 & 0 \\
rumba      & 1 & 0 & 0 & 1 & 0 & 0 & 0 & 1 & 0 & 0 & 1 & 0 & 1 & 0 & 0 & 0 \\
shiko      & 1 & 0 & 0 & 0 & 1 & 0 & 1 & 0 & 0 & 0 & 1 & 0 & 1 & 0 & 0 & 0 \\
son        & 1 & 0 & 0 & 1 & 0 & 0 & 1 & 0 & 0 & 0 & 1 & 0 & 1 & 0 & 0 & 0 \\
soukous    & 1 & 0 & 0 & 1 & 0 & 0 & 1 & 0 & 0 & 0 & 1 & 1 & 0 & 0 & 0 & 0
\end{tabular}
\label{tab:afro-cuban}
\end{center}
\end{table}

Considering the six rhythmic patterns of Table \ref{tab:afro-cuban} as 6 points in the space, the 
barcode is shown in Figure E1. 
The lack of Betti numbers for $H_i$ with $i > 0$ is indicative
of the underlying curiosity of this data as a metric space -- the distance between the
points is as if the points were the vertices of a graph with the distance being the  
{\it metric graph} distance.  Precisely, given a graph whose edges are labelled with 
positive numbers, the {\it metric distance} between two vertices is the minimum (over
all paths) of the sum of all edge lengths along a path from one point to the other. 
In this case it is a metric graph where all edges have length $0.0625$. 
See Figure E2.

\section{Conclusions and Discussion} \label{sec:conclusions}

There are two main contributions of this paper. First, it augments the current state of the art in understanding the topological structure of musical spaces by presenting several new and expanded examples of distance functions and the spaces they imply. The techniques are drawn from classical topology, and may be readily applied to a wide variety of distance measures. The development emphasizes the importance that even small details of the distance function can have on the global structure of the space, and shows that (generalized) circles and discs appear frequently among the spaces because of the dual nature of the distance functions: the pitch-based distances measure proximity in absolute pitch and proximity after octave reduction, the temporal distances measure proximity in time and proximity in phase around the rhythmic cycle. 

The second contribution is the application of persistent homology to musical data. Time-delay embeddings allow the method to exploit the time-dependent nature of musical data. Using this tool, the Betti numbers (and hence homology classes) can be calculated directly from time-pitch data. To verify that the tool makes sense, we have presented a number of cases where features are known and then verified that the topological analysis is able to recover those features using only the note-level data. As far as we know, this is the first time such a validation has been attempted. When looking at more complicated musical data, the results of the low dimension homology classes are usually easy to interpret.  In contrast, the meaning of the higher-dimensional Betti numbers is less clear.  In part this is due to the developing nature of persistent homology as a field, and we suspect that once the persistent homology of Gaussian fields is better understood, it will help in the understanding of barcodes such as those presented here. 

\bibliographystyle{IEEE}

\newpage
\section*{Electronic Supplement to {\em Topology of Musical Data}}

This document is the electronic supplement for the article 'Topology of Musical Data' 
by Sethares and Budney.

\begin{appendix}

\section{Example \ref{ex:s1xs2}}

\begin{figure}
\begin{center}{
\psfrag{x}[tl][tl][1][0]{$x$} 
\psfrag{y}[tl][tl][1][0]{$y$} 
\psfrag{z}[tl][tl][1][0]{$z$} 
\psfrag{i}[tl][tl][1][0]{$\vec i$} 
\psfrag{j}[tl][tl][1][0]{$\vec j$} 
\psfrag{k}[tl][tl][1][0]{$\vec k$} 
\psfrag{i-j}[tl][tl][1][0]{$\vec i - \vec j$} 
\psfrag{j-k}[tl][tl][1][0]{$\vec j - \vec k$} 
\psfrag{j-i}[tl][tl][1][0]{$\vec j - \vec i$} 
\psfrag{k-i}[tl][tl][1][0]{$\vec k - \vec i$} 
\psfrag{k-j}[tl][tl][1][0]{$\vec k - \vec j$} 
\psfrag{i-j}[tl][tl][1][0]{$\vec i - \vec j$} 
\includegraphics[width=6cm]{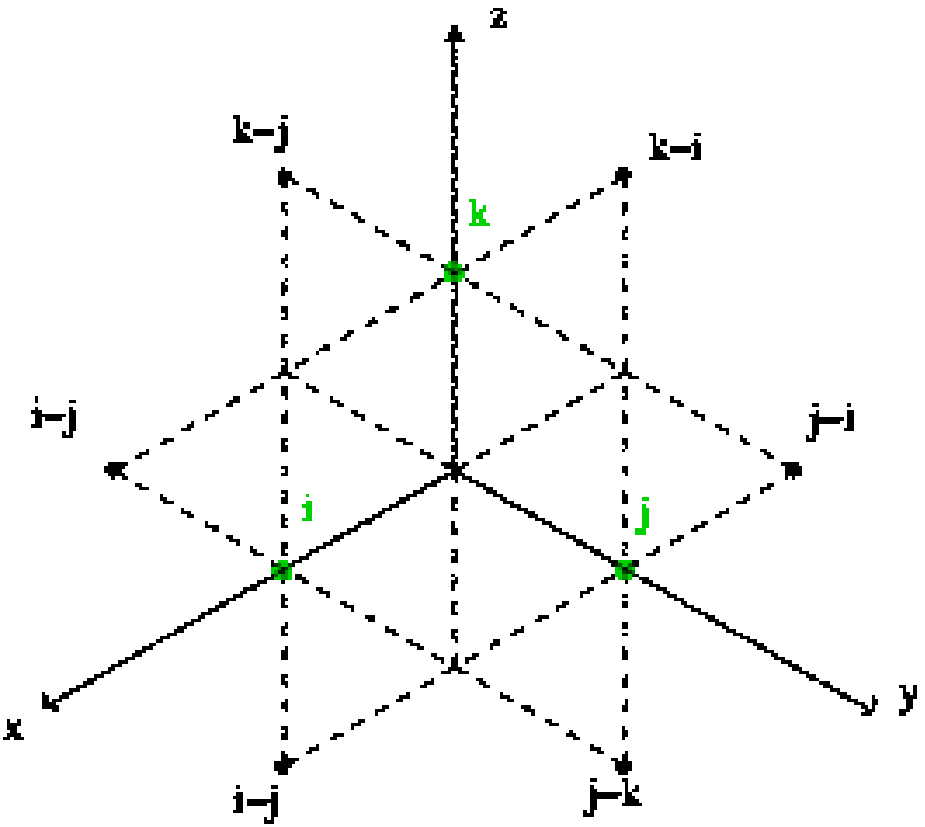}
\includegraphics[width=5cm]{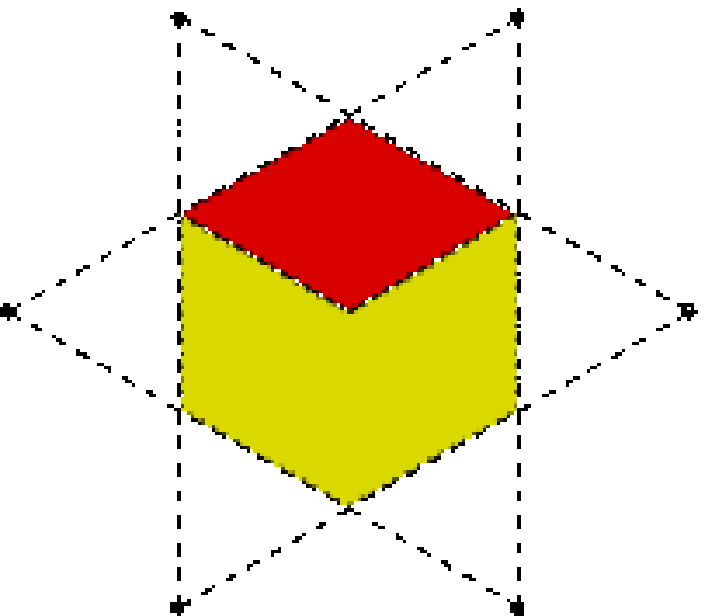}
}\end{center}
\caption{The fundamental domain for the 3-note melody space
is described in detail in the text.
\label{fig:fundamentalDomain}}
\end{figure}

For the subgroup $G=A_3 = \Zed_3$, $Symm^G_3(S^1)$ represents the 3-note melody space. 
Analogous to Morton's study of symmetric products, this is homeomorphic to 
$S^1 \times S^2$. This can be shown by recognizing that Morton's bundle 
(Theorem \ref{MortonThm}) $Symm^{A_3}_3(S^1) \to S^1$ is trivial, and the fibre is a sphere.  
To see this, observe that a fundamental domain in Figure \ref{fig:fundamentalDomain} 
for the action of $\Delta_\Zed \rtimes A_3$ on $\Delta_\Real$ is the quadrilateral with vertices
$(0,0,0)$, $(-2/3, 1/3, 1/3)$, $(1/3,-2/3,1/3)$ and $(-1/3,-1/3,2/3)$ (the red quadrilateral).  
The argument proceeds in two steps. 
First, $\Delta_\Zed \subset \Delta_\Real$ is the dual of the hexagonal lattice; 
it is the triangulation of the plane by equilateral triangles.  
Take the hexagon centered at the origin, dual to the triangular lattice. 
This has vertices $(-2/3, 1/3, 1/3)$, $(1/3,-2/3,1/3)$ together with their cyclic 
permutations (given by the action of $A_3$).   
Accordingly, the fundamental domain of the $\Delta_\Zed \rtimes A_3$ is one third of the hexagon, the quadrilateral described above.  Within the red polygon, the two top edges are 
glued together by a translation-rotation, while the two bottom edges are glued 
together by a counter-clockwise $4\pi/3$ rotation about the origin.  
From this perspective the fiber is naturally a triangulated $S^2$, 
using precisely two triangles.  The triangulation can be viewed as the 
union of two identical triangles along their common boundary.  
The monodromy has two fixed points (the centres of the triangles) and 
is rotation by $\pm 2\pi/3$ at those fixed points. 

Figure \ref{fig:fundamentalDomain} shows a view `down' the $\Delta_\Real$ plane in $\Real^3$, viewed from the perspective of the $x=y=z$ axis.  The green dots indicate the points lying outside the $\Delta_\Real$ plane, while the black dots are in the $\Delta_\Real$ plane.  The red quadrilateral is the fundamental domain of the $\Delta_\Zed \rtimes A_3$ action on $\Delta_\Real$, while the union of the yellow and red polyhedrons is the fundamental domain for the action of $\Delta_\Zed$ on $\Delta_\Real$ (the hexagonal lattice). 

A section of the bundle is realized by all the $3$-note melodies where all three notes equally spaced; these are the fixed-points of the monodromy $S^2 \to S^2$ above.  There are two such natural sections: one where the three notes are cyclically increasing, and one where the notes are cyclically decreasing in pitch.  Since the bundle is trivial, there are many other sections, but since the monodromy only has the two fixed points, no other sections are pitch-translates of a point in the fiber. 
This argument explains example \ref{ex:s1xs2}.

\section{An alternative representation} \label{poincare}

A common way of representing the Betti numbers of a space used in topology is called the Poincar\'e Polynomial.  Given a space $X$, the Poincar\'e Polynomial is the function $p_X(t) = \sum_{i=0}^\infty \beta_i(X) t^i$ where $\beta_i(X)$ is the $i$-th Betti number of $X$, i.e. $\beta_i(X) = rank(H_i X)$.  The Poincar\'e Polynomial of a circle $S^1$ is $p_{S^1}(t) = 1+ t$ since the only non-trivial homology groups are $H_0$ and $H_1$, both are free abelian on one generator.   Similarly, for the $n$-sphere $p_{S^n}(t) = 1+t^n$.  A consequence of the K\"unneth Theorem is that the Poincar\'e Polynomials of a product are the products of the Poincar\'e polynomials of the individual spaces. See \cite{Hatcher} for details of the K\"unneth Theorem and Poincar\'e polynomials.  
\[ p_{X \times Y} = p_X \cdot p_Y \]
Thus, in the case of Theorem \ref{prop:delayChords}, 
\[ p_{Symm^G_n S^1}(t) = (1+t)^j = \sum_{j=0}^n {n \choose j} t^j. \]
Similarly,
\[ p_{Symm^{A_3}_3 S^1}(t) = (1+t)(1+t^2) = 1+t+t^2+t^3. \]

\begin{example} \label{ex:Symm4ModZ4}
$Symm^{\Zed_4}_4 (S^1)$ can be analyzed via its Morton bundle.  The fibre of the map $Symm^{\Zed_4}_4(S^1) \to S^1$ is $\Delta_\Real / \Delta_\Zed \rtimes \Zed_4$.  Let $e_1,e_2,e_3,e_4$ represent the standard basis vectors in $\Real^4$.  The lattice $\Delta_\Zed \subset \Delta_\Real$ is generated by the set $\{e_i-e_j : i \neq j\}$, which makes the vertices of a cuboctahedra.  The tessellation of $\Delta_\Real$ dual to this lattice is the rhombic dodecahedral tessellation.  So a fundamental domain for the $\Delta_\Zed$-action on $\Delta_\Real$ is the rhombic dodecahedron centred at the origin, such that the centres of the faces correspond to the set $\{e_i-e_j : i \neq j\}$.  The action of $\Zed_4$ on the rhombic  dodecahedron is given by rotation by $\pi/2$ about an axis through the origin and a vertex of valence $4$, followed by the antipodal map.   The fundamental domain of the $\Delta_\Zed \rtimes \Zed_4$-action on $\Delta_\Real$ is a parallepiped.  This can be seen by expressing the rhombic dodecahedron as a union of $4$ isometric parallepipeds which $\Zed_4$ permutes cyclically.  The parallelepipeds correspond to the vertices of the regular tetrahedron vertices, the orthogonal projects of $\{e_1,e_2,e_3,e_4\}$, in particular the intersection of the four parallelpipeds consists of only the origin.   The identifications on the boundary of the parallelpiped folds two of the faces on themselves, and identifies two pairs of faces.  This leaves the fiber with a CW-structure with $4$ $0$-cells, $3$ $1$-cells, $2$ $3$-cells and a single $3$-cell.  The $2$-cell attaching maps are null, so the fiber has the homotopy-type of the wedge $S^2 \vee S^2$ attach a $3$-cell via the sum of twice the characteristic maps of the $2$-spheres respectively.  The result is the fibre $\Delta_\Real / \Delta_\Zed \rtimes \Zed_4$ has the homology groups $H_0 \cong \Zed$, $H_1 \cong 0$, $H_2 \cong \Zed \oplus \Zed_2$ and $H_3 \cong 0$. In particular, $Symm^{\Zed_4}_4(S^1)$ is not a manifold.   Since non-torsion $H_2$ elements are in the image of the quotient map $H_2 Symm_4(S^1) \to H_2 Symm^{\Zed_4}_4(S^1)$, the monodromy acts trivially on the homology of the fiber. So the homology groups of $Symm^{\Zed_4}_4(S^1)$ are isomorphic to that of a product bundle, $H_0 \cong \Zed$, $H_1 \cong \Zed$, $H_2 \cong \Zed \oplus \Zed_2$, $H_3 \cong \Zed \oplus \Zed_2$, and $H_4 \cong *$. 
\end{example}

\section{``Abbott's Bromley Horn Dance''} \label{sec:BromleysChords3D}

\begin{figure}
\centering \includegraphics[width=4.8in]{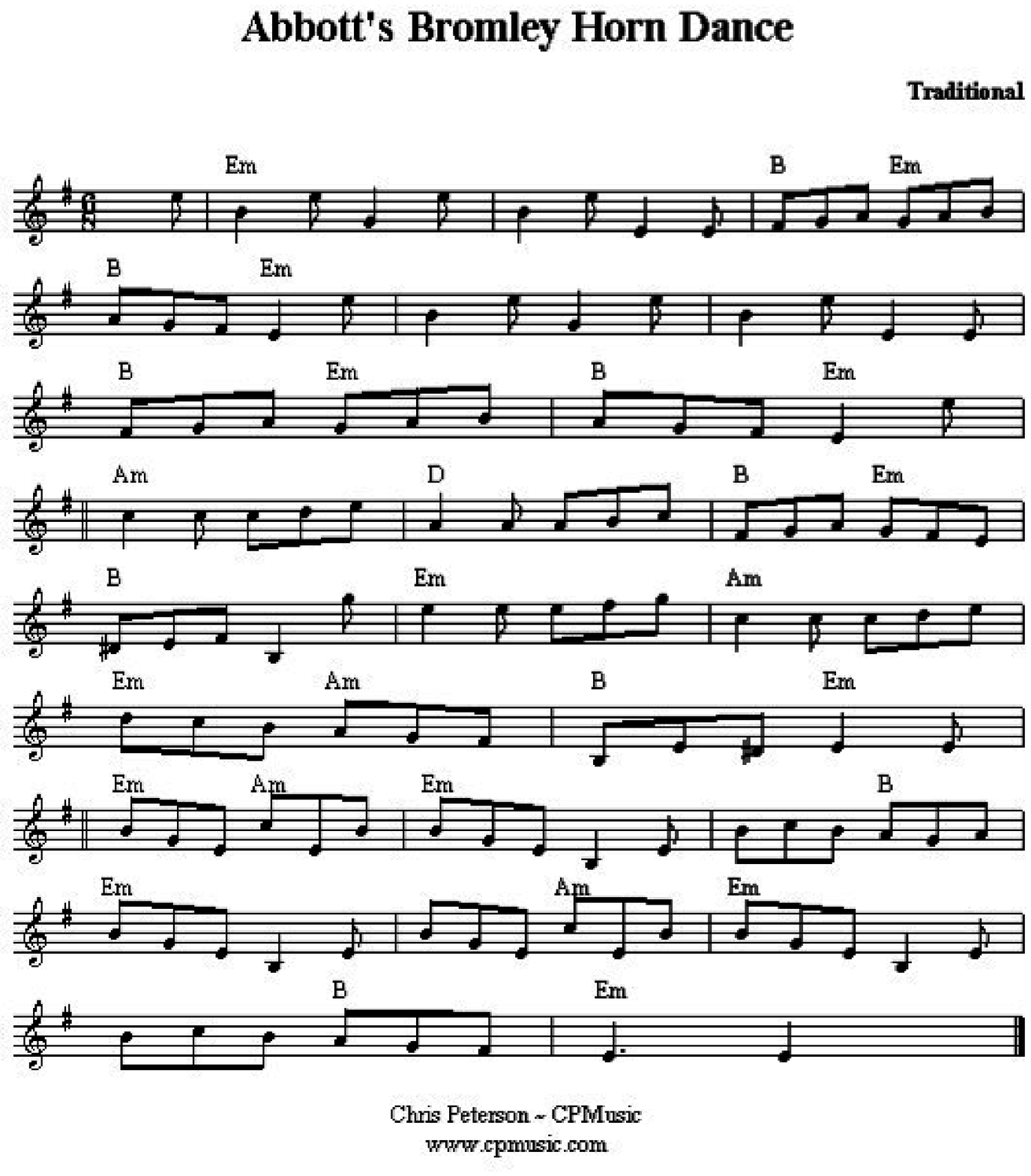}
\caption{The traditional melody ``Abbott's Bromley
Horn Dance'' is taken from Chris Peterson's collection \cite{peterson}. 
The standard MIDI version of this melody is analyzed using 
the ideas of persistent homology in Figures \ref{fig:abbotts1D},
\ref{fig:abbotts2D}, \ref{fig:abbotts3D},
and \ref{fig:BromleysChords3D}.  
\label{fig:abbottssheet}}
\end{figure}

\begin{figure}
\centering \includegraphics[width=4.8in]{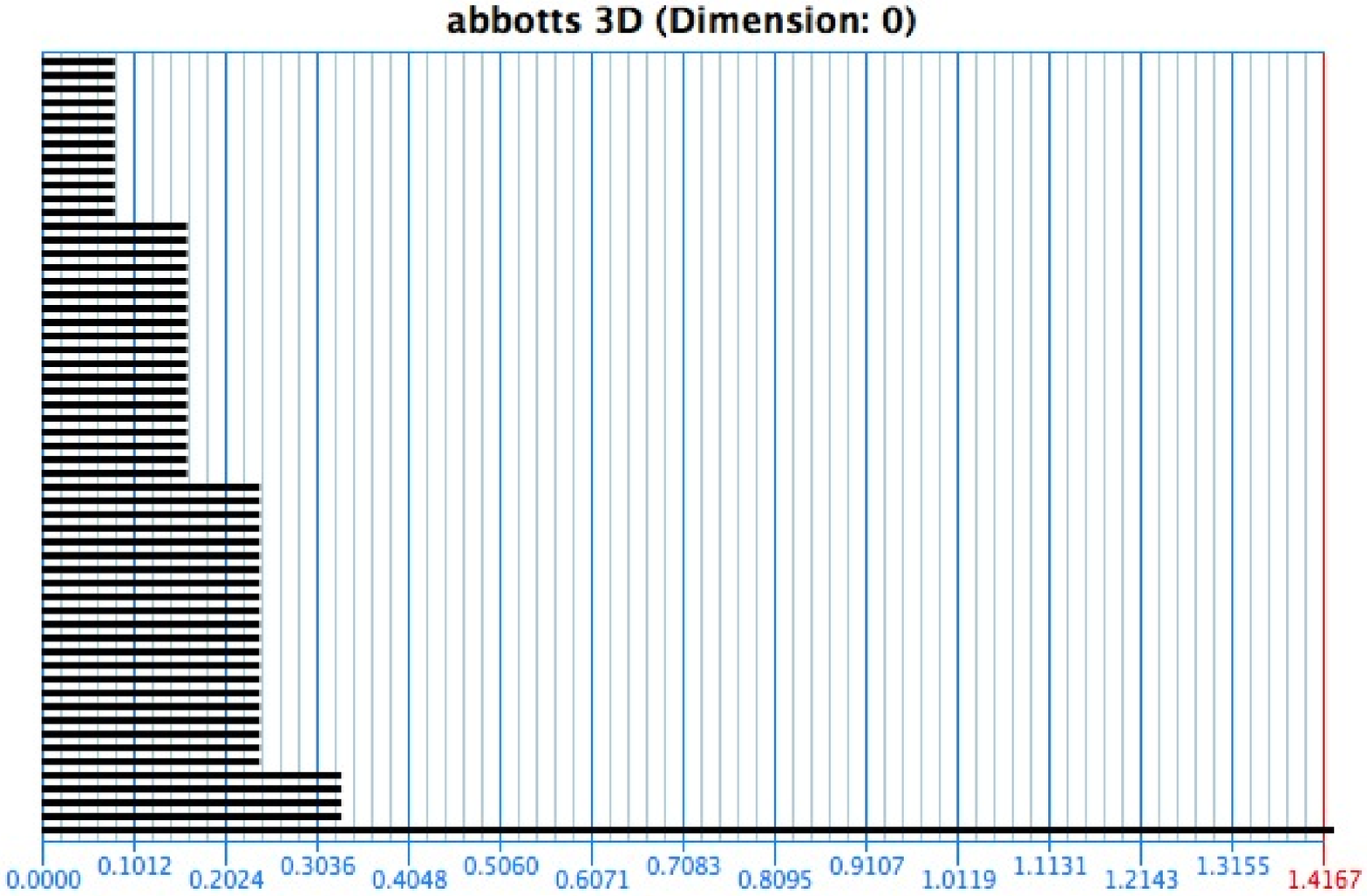}
\includegraphics[width=4.8in]{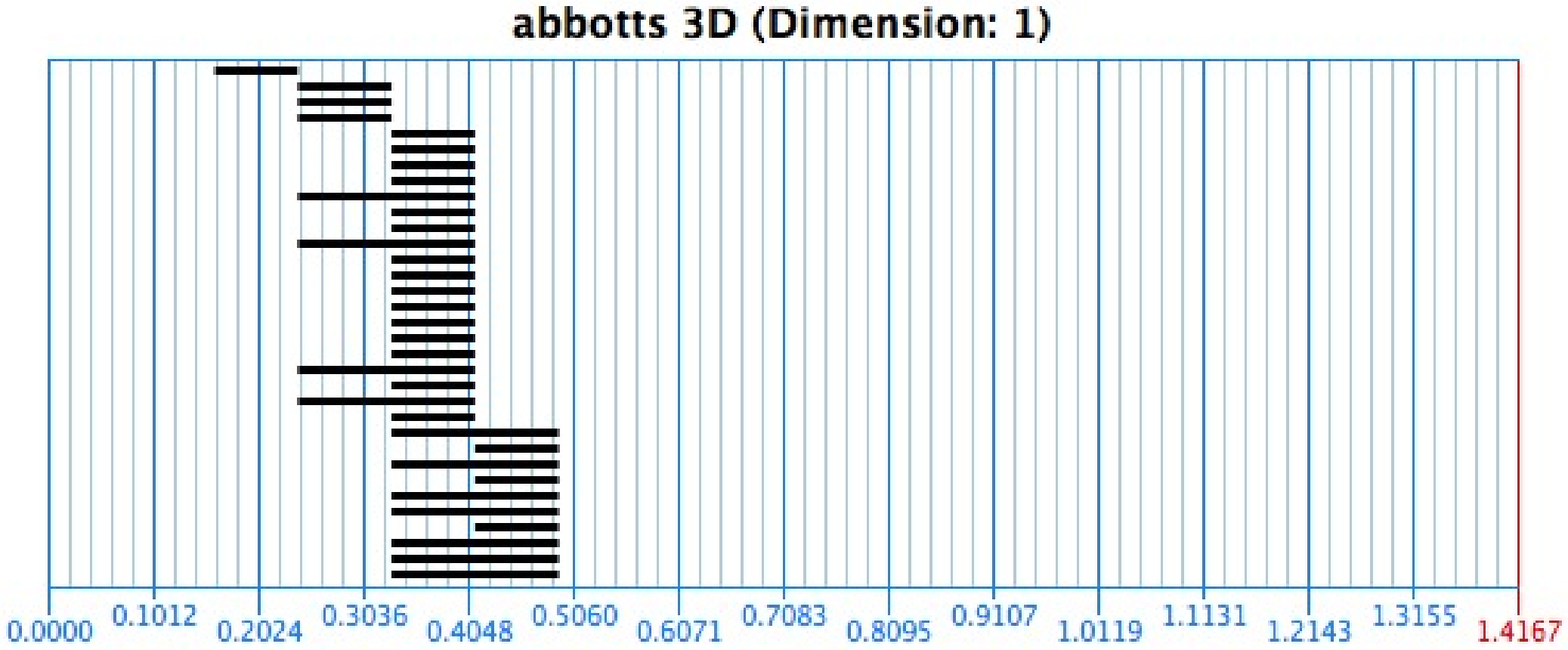}
\includegraphics[width=4.8in]{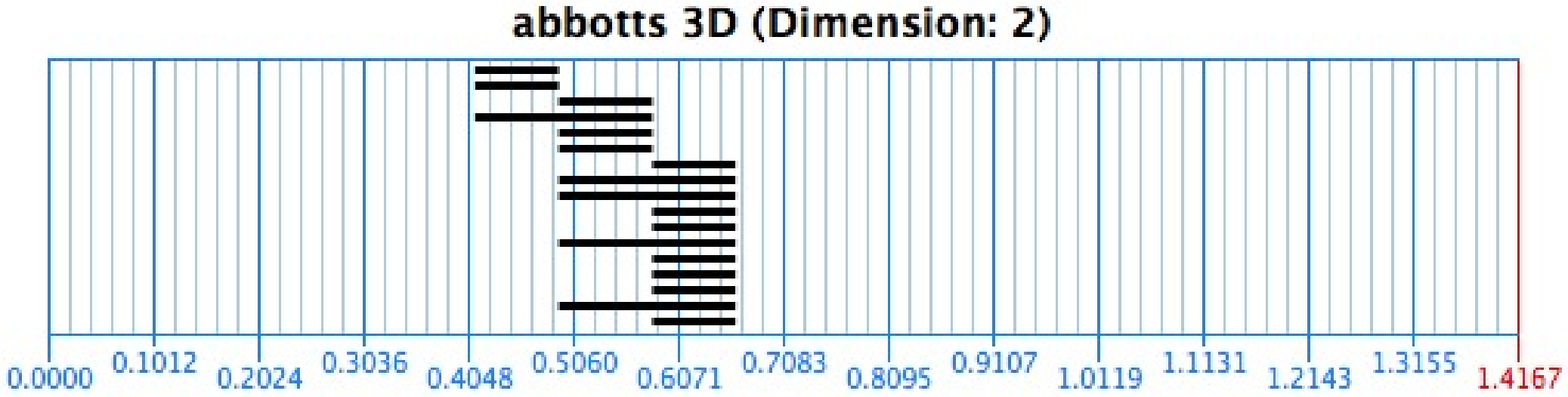}
\includegraphics[width=4.8in]{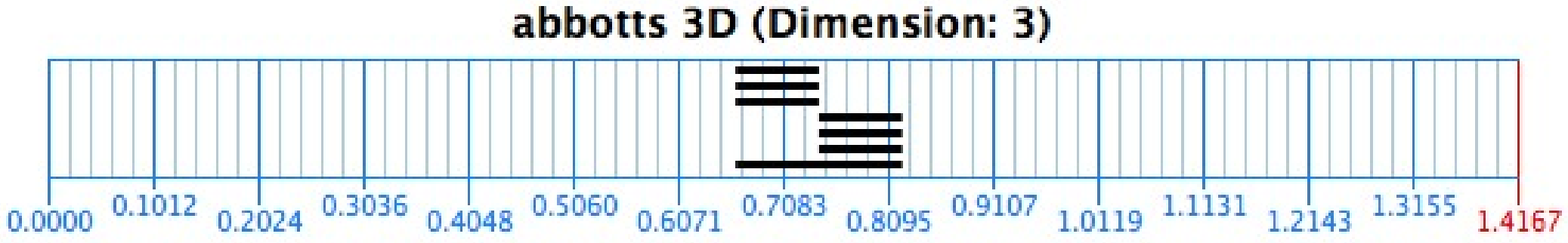}
\caption{Barcodes for the two-dimensional time-delay 
embedding of ``Abbott's Bromley Horn Dance'' (see Figure \ref{fig:abbottssheet}) 
show a remarkable array of features that represent both the structure of
the space and the structure of the individual piece. 
The dimension 0 plot is straightforward (showing
the number of melodic 3-tuples at each distance), the Dimension 1 and
2 plots contain a fascinating collection of circles and higher 
dimensional analogs that persist over a nontrivial range of $\epsilon$. 
The distance function is (\ref{eqn:distancePC}) applied to 
each element of the 3-vector, and then summed.
The topological space in
which this is embedded is $S^{1}\times S^{1}\times S^{1}$, the 
$n=3$ case of Example \ref{ex:orderedTuples}.
The maximum distance between elements of the performance is 
$1.4$, which is close to the diamater 1.5 of the complete space.
\label{fig:abbotts3D}}
\end{figure}

\begin{figure}
\centering \includegraphics[width=4.8in]{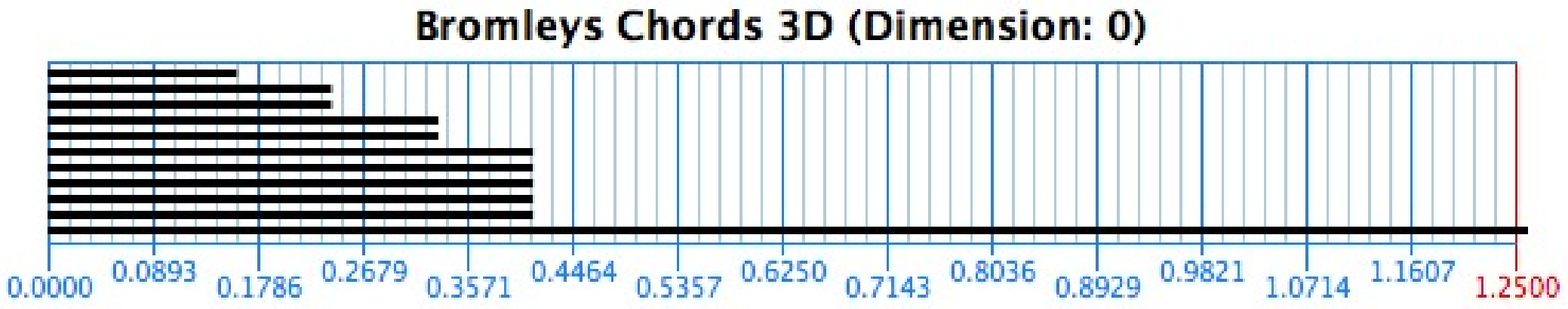}
\includegraphics[width=4.8in]{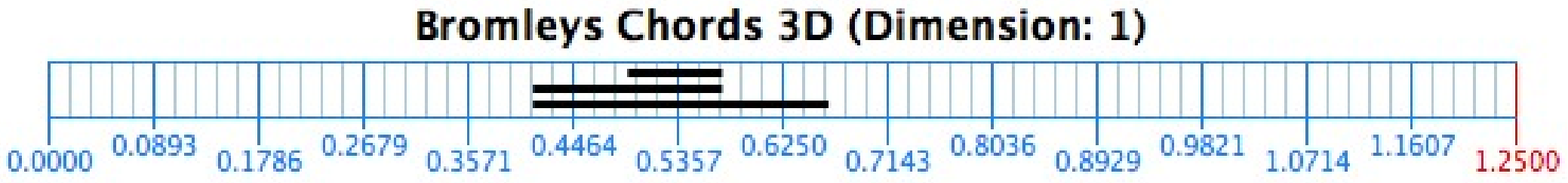}
\caption{Barcodes for the dimension-two time-delay
embedding of the chord progression 
from the score of ``Abbott's Bromley Horn Dance.''
The distance function is the two dimensional
time-delay embedding (\ref{eqn:distanceTDE}).
Under this distance, the diamater of the space is 4.5 and
the maximum distance between elements is 1.25, suggesting that
the data sits in a small `ball' and occupies only a small portion of
the full space. 
\label{fig:BromleysChords3D}}
\end{figure}

\section{Bach's Chorale No. 19} \label{sec:BachChorale19}

\begin{figure}
\centering \includegraphics[width=4.8in]{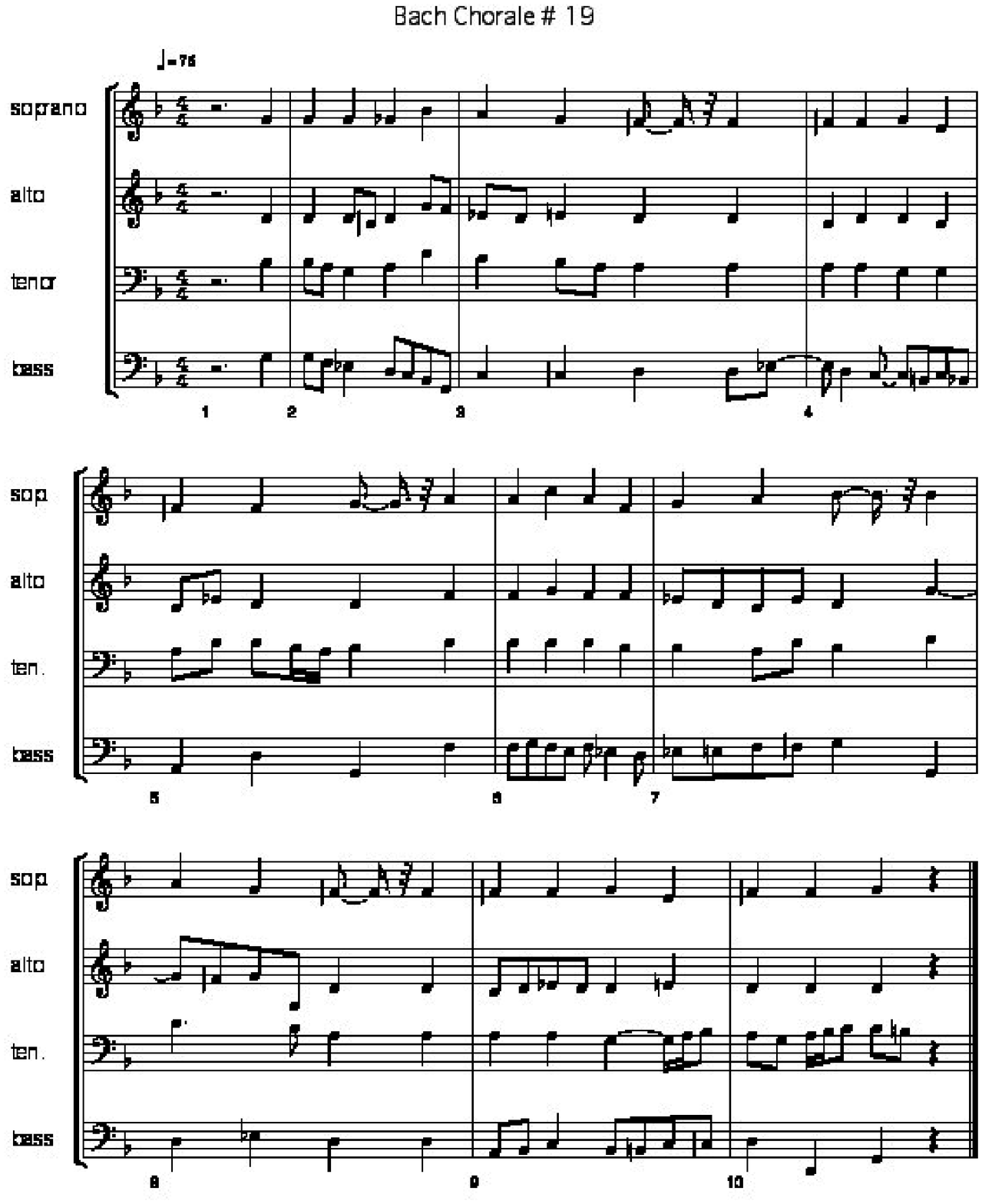}
\caption{The standard MIDI version of Bach's Chorale No. 19 is analyzed using 
the ideas of persistent homology in Figures \ref{fig:BachChorale19-1D}
and \ref{fig:BachChorale19-2D}.  
\label{fig:BachChorale19}}
\end{figure}

\begin{figure}
\centering \includegraphics[width=4.8in]{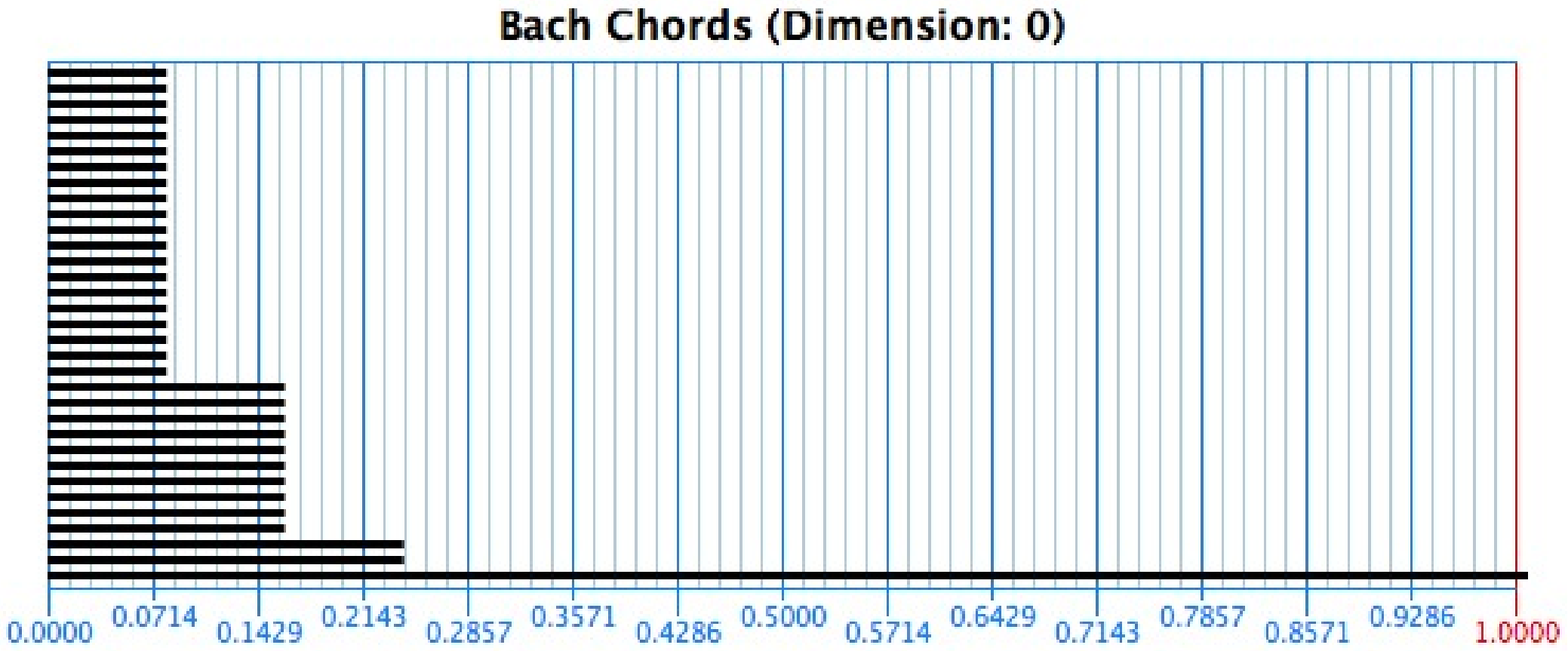}
\includegraphics[width=4.8in]{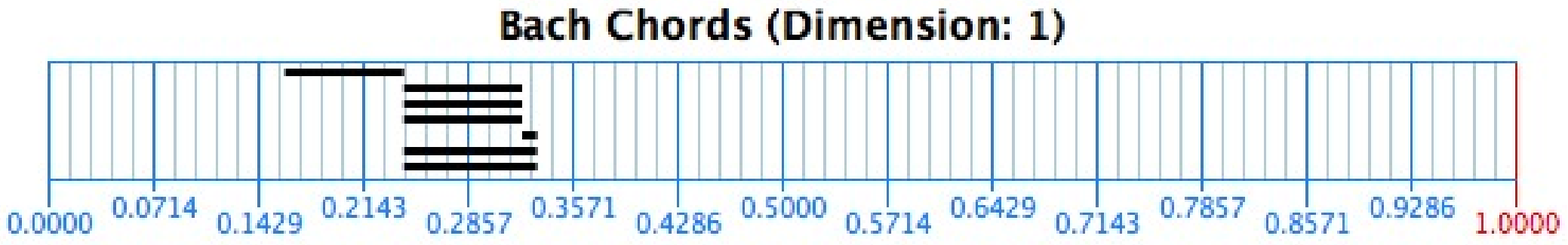}
\includegraphics[width=4.8in]{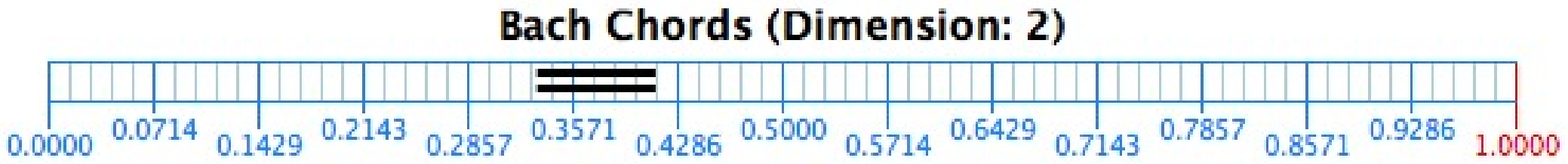}
\caption{Barcodes for Bach's Chorale No. 19
under the chord class distance (\ref{eqn:distanceCC}).
The space is the nonorientable bundle $S^{1} \times D^{3}$ 
(this is the $n=4$ case of Example \ref{ex:chordSpace})
and the total diamater of the space is 2.0 while  
the maximum distance between any two elements is 1.0.  This indicates
that the data sits in a small `ball' and fills only a small percent of
the full space of chords. 
\label{fig:BachChorale19-1D}}
\end{figure}

\begin{figure}
\centering \includegraphics[width=4.8in]{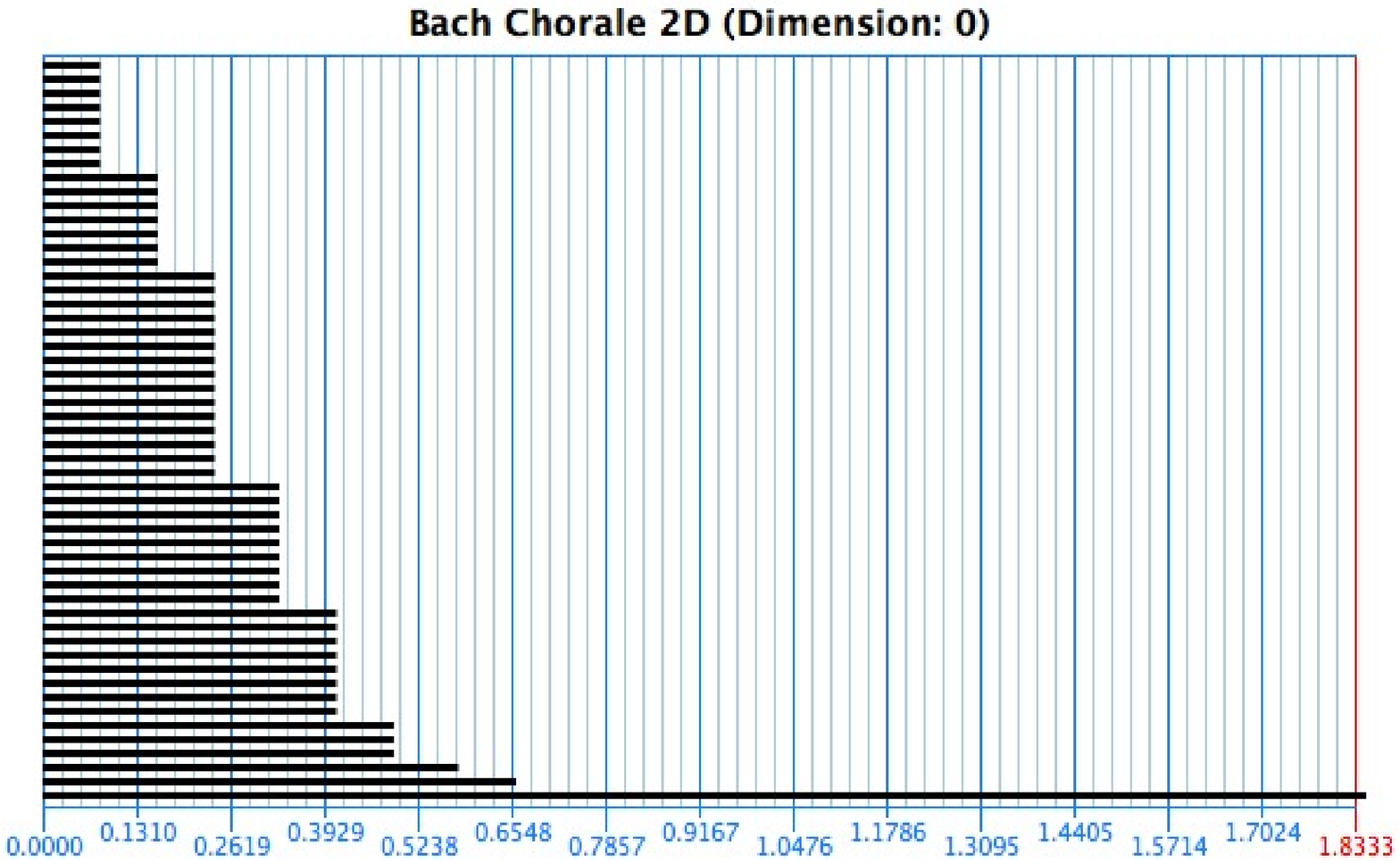}
\includegraphics[width=4.8in]{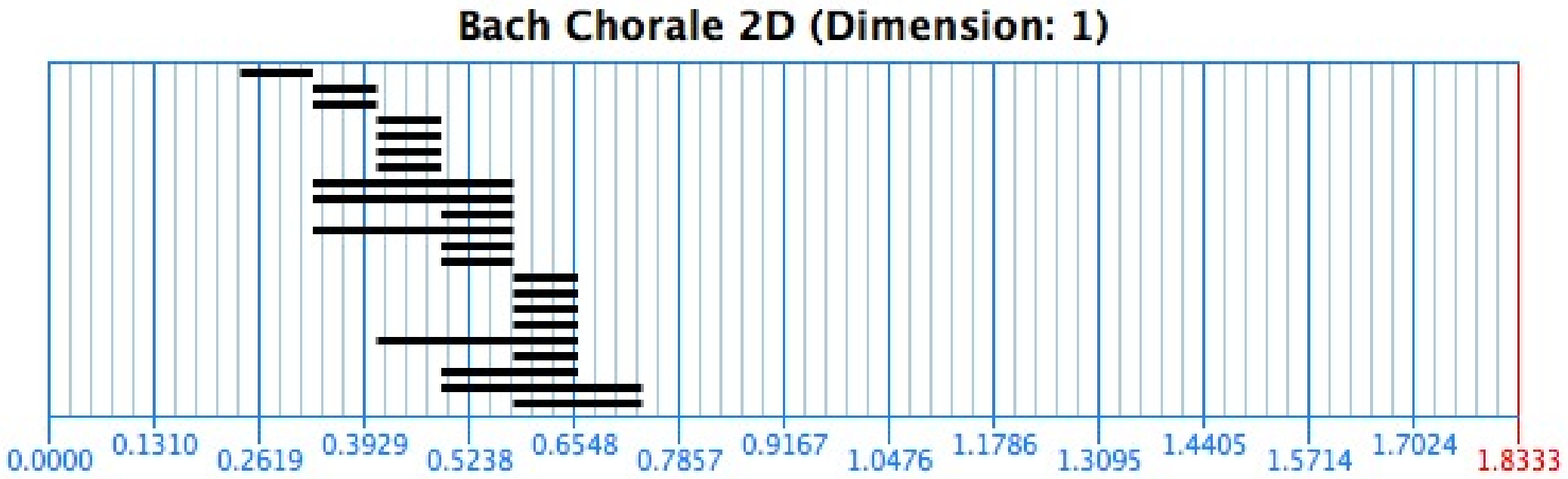}
\includegraphics[width=4.8in]{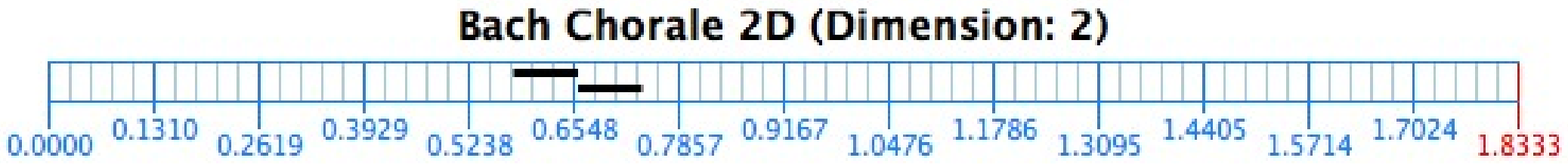}
\includegraphics[width=4.8in]{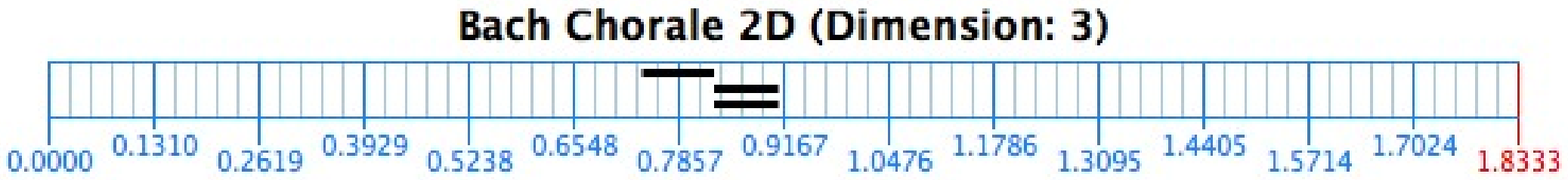}
\caption{Barcodes for the two-dimensional time-delay embedding
of Bach's Chorale No. 19.
The distance function is (\ref{eqn:distanceTDE}). According to Proposition
\ref{prop:delayChords} with 
time delay $j=2$ and chord-size $k=4$.
The diameter of the space is $4$ and the maximum distance between
elements is $1.833$. 
\label{fig:BachChorale19-2D}}
\end{figure}

\section{Rhythmic barcodes} \label{sec:afro-cuban}

\begin{figure}
\centering \includegraphics[width=4.8in]{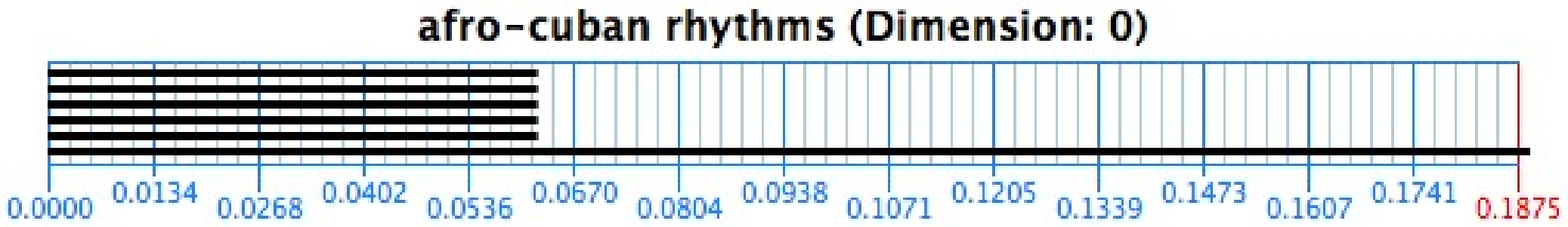}
\caption{Barcodes for the six rhythmic patterns in Table \ref{tab:afro-cuban} under the
rhythm distance (\ref{eqn:distanceCyc}). The max distance between elements is
$0.1875$. 
\label{fig:afro-cubanrhythms0}}
\end{figure}

\begin{figure}
{
\psfrag{0}[tl][tl][1][0]{BN}
\psfrag{1}[tl][tl][1][0]{G}
\psfrag{2}[tl][tl][1][0]{R}
\psfrag{3}[tl][tl][1][0]{Sh}
\psfrag{4}[tl][tl][1][0]{CS}
\psfrag{5}[tl][tl][1][0]{So}
$$\includegraphics[width=2in]{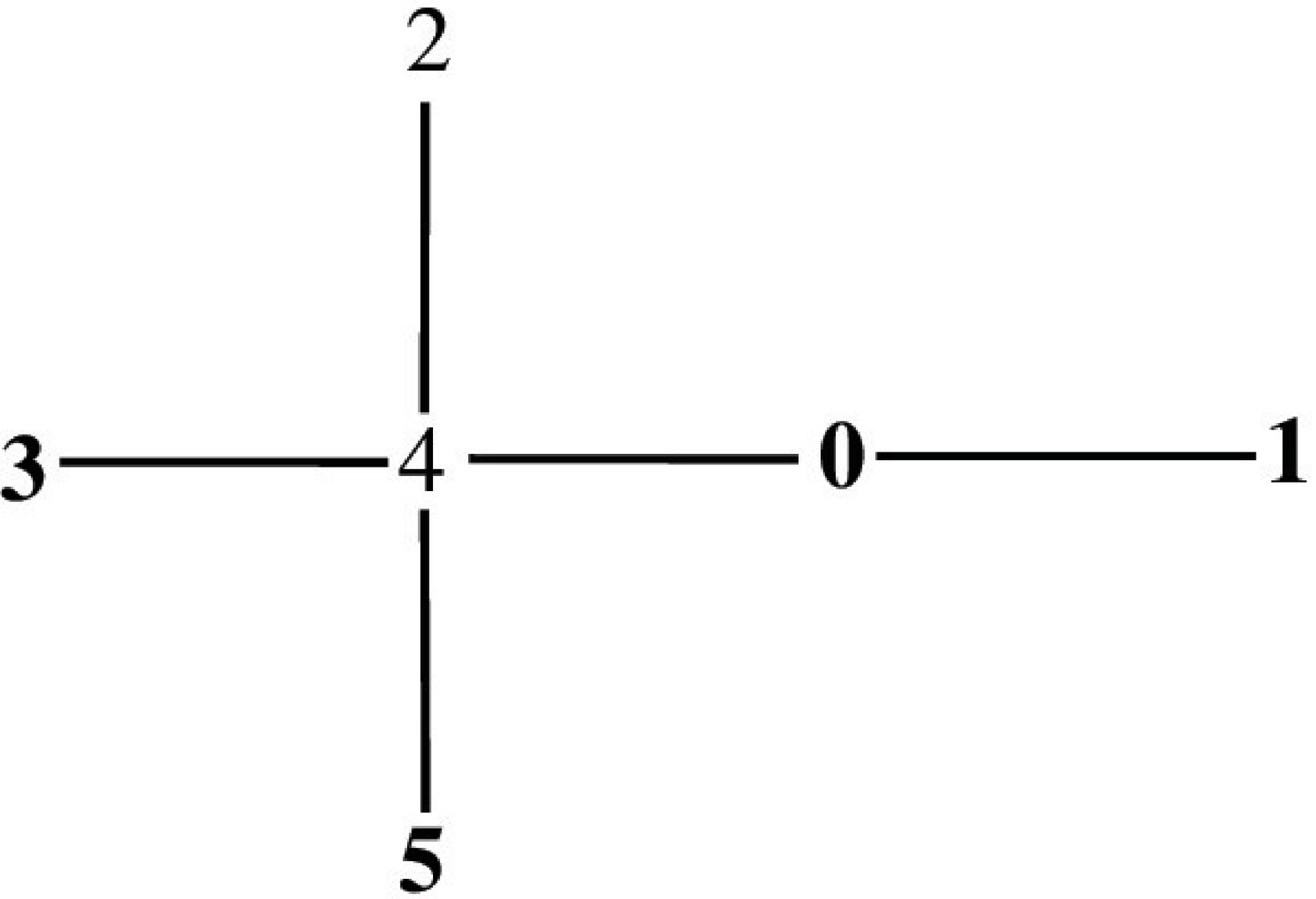}$$
}
\caption{Representation of the underlying metric graph describing the data.}
\label{fig:afro-cubanrhythms1}
\end{figure}

\section{Topology of Musical Data}

\end{appendix}
\end{document}